\documentclass[11pt]{article}
\usepackage[top = 1.2in, bottom= 1.2in, left = 1.3in, right=1.3in]{geometry}

\usepackage{paper1-macros}
\title{A Geometric transition from hyperbolic to anti de Sitter geometry}
\author{Jeffrey Danciger\footnote{J.D. is partially supported by the National Science Foundation, grant number DMS 1103939}\\
  Department of Mathematics,\\
  University of Texas - Austin,\\
  1 University Station C1200,\\
  Austin, TX 78712-0257\\
  \texttt{jdanciger@math.utexas.edu}}
\date{\today}

\begin{document}

\removenote

\maketitle

\begin{abstract}

We introduce a geometric transition between two homogeneous three-dimensional geometries: hyperbolic geometry and anti de Sitter (AdS) geometry. Given a path of three-dimensional hyperbolic structures that collapse down onto a hyperbolic plane, we describe a method for constructing a natural continuation of this path into AdS structures. In particular, when hyperbolic cone manifolds collapse, the AdS manifolds generated on the ``other side" of the transition have tachyon singularities. The method involves the study of a new transitional geometry called \emph{half-pipe} geometry. We demonstrate these methods in the case when the manifold is the unit tangent bundle of the $(2,m,m)$ triangle orbifold for $m \geq 5$.
\end{abstract}

\section{Introduction}

A geometric structure on a manifold is a system of local coordinates modeled on a homogeneous space $X = G/K$. The structure is called \emph{hyperbolic}, if modeled on hyperbolic space $\HH^n = \PO(n,1)/\OO(n)$, or \emph{AdS} if modeled on the anti de Sitter space $\AdS^n = \PO(n-1,2)/\OO(n-1,1)$, which is a Lorentzian analogue of $\HH^n$.
We study families of geometric structures that \emph{degenerate}, or collapse, our main source of examples coming from three-dimensional hyperbolic cone manifolds. The guiding philosophy, based loosely on Thurston's geometrization program, is that when deforming structures collapse, the degeneration is signaling that a transition to a different geometry is needed in order to continue the deformation; the form of the collapse should tell us which type of geometry lies on the other side of the transition. This article develops a new geometric transition,  going from hyperbolic structures 
to $\AdS$ structures, in the context of structures that collapse onto a co-dimension one hyperbolic plane.

Let $N$ be a closed three-manifold and let $\Sigma \subset N$ be a knot. Assume that $N$ itself is not hyperbolic, but that $N \setminus \Sigma$ admits a complete hyperbolic structure.
Then there exist hyperbolic cone structures on $N$ with cone singularity along $\Sigma$. These are parametrized by the cone angle $\alpha$ and in many cases all cone angles $\alpha < 2\pi$ may be achieved.
As $\alpha \to 2\pi^-$, however, these structures must degenerate and in general the collapse may take a variety of forms. 
If, for example, the manifolds collapse uniformly to a point, then a transition to spherical geometry passing through Euclidean may be possible.
This well-known transition, studied by Hodgson \cite{Hodgson-86}, Porti \cite{Porti-98}, and later in the proof of the Orbifold Theorem by Cooper-Hodgson-Kerckhoff~\cite{Cooper-00} and Boileau-Leeb-Porti~\cite{Boileau-05}, can be described as follows.
The collapsing hyperbolic metrics are rescaled to give a limiting Euclidean metric and that Euclidean metric is then deformed to give nearby spherical cone metrics.
Porti and collaborators have also studied degenerations in the context of non-uniform collapse to a point \cite{Porti-02}, and collapse to a line \cite{Huesener-01}.

Consider a path of hyperbolic cone structures on $N$ such that each chart in $\HH^3$ is collapsing onto a particular copy $\plane$ of $\HH^2$. 
The transition functions between charts, which lie in $\text{Isom} (\mathbb H^3) = \PO(3,1)$, are converging into the $\OO(2,1)$ subgroup that preserves $\plane$.  
The collapsed charts define a transversely hyperbolic foliation $\mathcal F$ on $N$. Collapse of this type may happen, for example, when $N$ is Seifert fibered over a base $S$ of hyperbolic type; in this case, one expects collapse onto a hyperbolic representative of the surface $S$ so that $\mathcal F$ agrees with the Seifert fibration.  
The problem of \emph{regeneration}, or recovery, of nearly collapsed hyperbolic structures from a transversely hyperbolic foliation was examined by Hodgson~\cite{Hodgson-86}, and later in a specific case by Porti~\cite{Porti-10}. However, it had not yet been established how to construct a geometric transition in this context.
Our point of view, based on real projective geometry, is that such degeneration behavior suggests a natural transition to anti de Sitter geometry. 
As motivation, note that both hyperbolic and $\AdS$ are sub-geometries of projective geometry and that, when suitably embedded,
the structure groups $\Isom \HH^3 = \PO(3,1)$ and $\Isom \AdS^3 = \PO(2,2)$ intersect inside $\PGL(4,\RR)$ in the $\OO(2,1)$ subgroup that preserves a hyperbolic plane $\plane$.
We develop tools to construct examples of this hyperbolic-$\AdS$ transition. In the setting above, in which the hyperbolic structures have cone singularities, the $\AdS$ structures generated on the other side of the transition will have the analogous singularity in Lorentzian geometry, called a tachyon. We demonstrate our methods in the case that $N$ is the unit tangent bundle of the $2,m,m$ triangle orbifold.

\begin{Theorem}
\label{thm:regen-2mm}
Let $m \geq 5$, and let $S$ be the hyperbolic structure on the two-sphere with three cone points of order $2,m,m$.  Let $N$ be the unit tangent bundle of $S$. Then, there exists a knot $\Sigma$ in $N$ and a continuous path of real projective structures $\mathscr P_t$ on $N$, singular along $\Sigma$, such that $\mathscr P_t$ is conjugate to:
\begin{itemize}
\item
a hyperbolic cone structure of cone angle $\alpha < 2\pi$, when $t > 0$, or
\item
an $\AdS$ structure with tachyon singularity of mass $\phi < 0$, when $t < 0$.
\end{itemize}
As $t \to 0$, the cone angle $\alpha \to 2 \pi^-$ (resp. $\phi \to 0^-$) and the hyperbolic geometry (resp. $\AdS$ geometry) of $\mathscr P_t$ collapses to the surface $S$. 
\end{Theorem}

The key to the proof of Theorem~\ref{thm:regen-2mm}, given in Section~\ref{example1}, is the construction of the \emph{transitional} projective structure $\mathscr P_0$. We build $\mathscr P_0$ out of a new homogeneous geometry called \emph{half-pipe} ($\HP$) geometry, designed to bridge the gap between hyperbolic and $\AdS$. Half-pipe geometry arises naturally as a limit of conjugates of hyperbolic geometry inside of projective geometry. In the context above, in which hyperbolic structures collapse onto a hyperbolic plane $\plane$, we \emph{rescale} $\HH^3$
in order to prevent collapse. Specifically, we apply a projective transformation that preserves $\plane$ (the plane of collapse) but stretches transverse to $\plane$. The limit of the conjugate copies of $\HH^3$ produced by this rescaling is the projective model for half-pipe geometry $\HP^3$. We note that information about the collapsed hyperbolic structure is contained in the $\HP$ geometry: Every $\HP$ structure has a natural projection to its associated transversely hyperbolic foliation $\mathcal F$. Additionally, the $\HP$ geometry preserves some first order information in the collapsing direction; this information is critical for regeneration to both hyperbolic and $\AdS$ structures.

We emphasize that the most important contribution of this paper is the general construction of half-pipe geometry, given in Section~\ref{HP}, as the transitional geometry connecting hyperbolic to $\AdS$.
Although our main application is to structures with cone singularities, this construction is very general and may be useful in other contexts. We also note that half-pipe geometry does not have an invariant metric (Riemannian nor Lorentzian). Our projective geometry approach eliminates metrics from the analysis entirely.

\subsection{Cone/tachyon transitions}

Recall that a cone singularity (see Section~\ref{hyp-cone} or e.g. \cite{Hodgson-05, Bromberg-07, Suoto-03}) in $\HH^3$ geometry is a singularity along a geodesic axis such that the holonomy of a meridian encircling the axis is a rotation around the axis. 
Similarly, a tachyon (see Section~\ref{AdS-tachyon} or  \cite{Barbot-09}) in $\AdS^3$ geometry is a singularity along a space-like axis such that the holonomy of a meridian encircling the axis is a Lorentz boost orthogonal to the axis. The magnitude of the boost is called the tachyon mass. When rescaled, a hyperbolic cone singularity of cone angle approaching $2 \pi$ and an $\AdS$ tachyon singularity of mass approaching zero both limit to an \emph{infinitesimal cone singularity} in $\HP$. 

   As is the case for the hyperbolic-spherical transition, it is technically difficult to prove general statements about when a collapsing family of hyperbolic cone manifolds yields a transition to $\AdS$; we do not prove such statements here.
Rather we shift our point of view toward constructing geometric transitions from an $\HP$ structure (regeneration rather than degeneration). One main tool for such constructions, used in particular to prove Theorem~\ref{thm:regen-2mm}, is the following regeneration theorem.

\begin{Theorem}
\label{thm:main-regen}
Let $N$ be a closed orientable three-manifold with a half-pipe structure $\mathscr H$, singular along a knot $\Sigma$ and with infinitesimal cone angle $-\omega < 0$. 
Let $\rho_0: \pi_1 M \to \OO(2,1)$ be the holonomy representation of the associated transversely hyperbolic foliation $\mathcal F$. Suppose that $H^1(\pi_1 N \setminus \Sigma, \so(2,1)_{Ad \rho_0}) = \mathbb R$. Then there exists geometric structures on $N$, singular along $\Sigma$, and parametrized by $t \in (-\delta, \delta)$ which are
\begin{itemize}
\item hyperbolic cone structures with cone angle $2\pi - \omega t$ for $t > 0$
\item $\AdS$ structures with a tachyon of mass $\omega t$ for $t < 0$.
\end{itemize}
Both families limit, as projective structures, to the initial $\HP$ structure $\mathscr H$ as $t\to 0$.
\end{Theorem}

The proof of Theorem~\ref{thm:main-regen}, given in Section~\ref{singular}, involves a generalization of cone singularities to projective geometry. These \emph{cone-like} singularities include cone singularities in hyperbolic geometry, tachyons in $\AdS$ geometry, and infinitesimal cone singularities in $\HP$ geometry. 
The basic ingredient needed for Theorem~\ref{thm:main-regen} is an open-ness principle (Proposition~\ref{deform-cone-like}) for projective structures with cone-like singularities akin to the Thurston-Ehresmann principle for non-singular structures. We mention that in contrast to the convex projective structures  
studied by Goldman~\cite{Goldman-90}, Choi--Goldman~\cite{Schoi-97},  Benoist~\cite{Benoist-04}, Cooper--Long--Tillman~\cite{Cooper-12}, Crampon--Marquis~\cite{Crampon-12} and others,
the projective structures appearing in this article do not arise as quotients of domains in $\RP^3$.

%
%
%

The assumptions in Theorem~\ref{thm:main-regen} are satisfied by a variety of examples, including examples coming from small Seifert fiber spaces (as in Theorem~\ref{thm:regen-2mm}) and Anosov torus bundles. We refer the reader to \cite{Danciger-11} for a careful analysis of the torus bundle case using ideal tetrahedra.
The cohomology condition, reminiscent of a similar condition appearing in Porti's regeneration theorem for Euclidean cone structures \cite{Porti-98}, is simply a way to guarantee smoothness of the $\OO(2,1)$ representation variety. Indeed, given an $\HP$ structure, our construction of a geometric transition really only requires that a transition exists \emph{on the level of representations} which is implied by (but does not require) smoothness. 
In Section~\ref{s:brings}, we study certain structures collapsing to a punctured torus for which the $\OO(2,1)$ representation variety is not smooth. Though we can still produce transitions, we also observe an interesting flexibility phenomenon in this case: A transitional $\HP$ structure can be deformed so that it no longer regenerates to hyperbolic structures.

Theorem~\ref{thm:regen-2mm} gives a limited class of manifolds on which transitions may be constructed. Many more examples can be constructed using the techniques of \cite{Danciger-11}\marginnote{update reference, once triangulation paper is done}, in which transitioning structures are built out of ideal tetrahedra. 

\subsection{Some context}
From the point of view of three-manifold topology, our construction is a bit surprising. 
For in light of the work of Porti and others (see above), one might expect that collapsing hyperbolic cone manifolds should limit to the appropriate Thurston geometry of the underlying manifold, which in the context of Theorem~\ref{thm:regen-2mm} is the Riemannian geometry of $\widetilde{\SL(2,\RR)}$. Instead, Theorem~\ref{thm:regen-2mm} produces a half-pipe structure as the limit. We note that while half-pipe geometry and $\widetilde{\SL(2,\RR)}$ geometry are incompatible, our results do not rule out the possibility that $\widetilde{\SL(2,\RR)}$ structures could be constructed as limits of hyperbolic structures in some other way. However, as of yet no such construction exists.

On the other hand, from the point of view of $\AdS$ geometry, our results are not surprising at all. Many parallels in the studies of hyperbolic and $\AdS$ geometry have appeared in recent years, beginning with Mess's classification of maximal $\AdS$ space-times \cite{Mess-07, Andersson-07} and its remarkable similarity to the Simultaneous Uniformization Theorem of Bers \cite{Bers-60} for quasi-Fuchsian hyperbolic structures. Also noteworthy is the Wick rotation/rescaling theory of Benedetti-Bonsante \cite{Benedetti-09} which gives a correlation between the boundaries of convex cores of constant curvature space-times and those of geometrically finite hyperbolic three-manifolds. 
Stemming from Mess's work, results and questions in hyperbolic and $\AdS$ geometry (see \cite{BBDGGKKSZ-12}) have begun to appear in tandem, suggesting the existence of a deeper link. 
Our geometric transition construction provides a concrete connection, strengthening the analogy between the two geometries.

\bigskip
\noindent{\bf Organization of the paper.} Section~\ref{s:preliminaries} covers preliminary material about the various geometric structures appearing in this work. Section~\ref{HP} gives the construction of half-pipe geometry. Section~\ref{singular} develops the theory of projective structures with cone-like singularities leading to the proof of Theorem~\ref{thm:main-regen}. Section~\ref{example1} gives the proof of Theorem~\ref{thm:regen-2mm}, and Section~\ref{s:brings} discusses an interesting flexibility phenomenon for transitional $\HP$ structures.

\bigskip
\noindent {\bf Acknowledgements.} I am grateful to Steven Kerckhoff for advising some of this work during my doctoral studies at Stanford University. I also thank Joan Porti, for suggesting the possibility of the phenomenon in Section~\ref{s:brings}, and Suhyoung Choi for several helpful discussions about the content of Section~\ref{example1}. Many helpful discussions occurred at the \emph{Workshop on Geometry, Topology, and Dynamics of Character Varieties} at the National University of Singapore in July 2010, and also at the special trimester \emph{Geometry and Analysis of Surface Group Representations} at  Institut Henri Poincar\' e in Paris, 2012. I also acknowledge support from the National Science Foundation grants DMS 1107452, 1107263, 1107367 ``RNMS: GEometric structures And Representation varieties" (the GEAR Network).

%
%

\section{Preliminaries}
\label{s:preliminaries}

In this section we give a description of the $(G,X)$ formalism that is ubiquitous in the study of locally homogeneous geometric structures. We will then describe real projective, hyperbolic, and $\AdS$ geometry using this formalism. In the final subsection, we describe transversely hyperbolic foliations.  

\subsection{Deforming $(G,X)$ structures}
Let $G$ be a group of analytic diffeomorphisms of $X$. Recall that a $(G,X)$ structure on a manifold $M$ is given by a collection of charts into $X$ with transition maps (being the restrictions of elements) in $G$. This data is equivalent to the data of a local diffeomorphism $D: \widetilde M \to X$, called a \emph{developing map}, which is equivariant with respect to the \emph{holonomy representation} $\rho: \pi_1 M \to G$ (\cite{Ehresmann-36} or see \cite{Thurston, Goldman-10}).

A smooth family of $(G,X)$ structures on a manifold $M$ with boundary is given by a smooth family of developing maps $D_t : \widetilde M \rightarrow X$ equivariant with respect to a smooth family of holonomy representations $\rho_t : \pi_1 M \rightarrow G$. 
Two deformations $D_t$ and $F_t$ of a given structure $D_0$ are considered equivalent if there exists a path $g_t \in G$ and a path $\Phi_t$ of diffeomorphisms \emph{defined on all but a small neighborhood of $\partial M$} so that $$D_t = g_t \circ F_t \circ \widetilde \Phi_t$$
where $\widetilde \Phi_t$ is a lift of $\Phi_t$ to $\widetilde M$ and we assume $g_0 = 1$ and $\widetilde \Phi_0 = \text{Id}$.
A \emph{trivial} deformation of $D_0$ is of the form $D_t = g_t \circ D_0 \circ  \widetilde \Phi_t$. In this case, the holonomy representations differ by a path of conjugations:
\begin{equation*}
\rho_t = g_t \rho_0 g_t^{-1}.
\end{equation*}
Such a deformation of the holonomy representation is also called trivial.
Let $\mathscr R(\pi_1 M; G)$ be the space of representations up to conjugation (we only consider points at which this quotient is reasonable). Let $\mathscr D(M; G,X)$ be the space of all $(G,X)$ structures on $M$ up to the equivalence described above. The following fact, known as the \emph{Ehresmann-Thurston principle}, is crucial for the study of deformations of $(G,X)$ structures. 
\begin{Proposition}[Ehresmann, Thurston]
\label{prop:deform-withboundary}
The map $\operatorname{hol}: \mathscr D(M; G,X) \rightarrow \mathscr R(\pi_1 M; G)$, which maps a $(G,X)$ structure to its holonomy representation, is (well-defined and) a local homeomorphism.
\end{Proposition}

We emphasize that the definition of $\mathscr D(M; G, X)$ above does not consider behavior at the boundary. In particular, given a structure $D_0$ with special geometric features at $\partial M$, Proposition~\ref{prop:deform-withboundary} may produce nearby $(G,X)$ structures with very different boundary geometry. Often (and it will always be the case in this paper), it is desirable to deform $(G,X)$ structures with control over the geometry at the boundary. This is the case, for example, in the study of hyperbolic cone structures (see \cite{Hodgson-98}). The Proposition is not strong enough in these cases, and a thorough study of the boundary geometry is needed. In Section~\ref{singular} we will pay careful attention to this issue as we deform from hyperbolic cone structures to $\AdS$ tachyon structures.

\subsection{Infinitesimal Deformations}
Consider a smooth family of representations $\rho_t: \pi_1 M \rightarrow G$. The derivative of the homomorphism condition (evaluated at $t=0$) gives that $$\rho'(ab) = \rho'(a)\rho_0(b) + \rho_0(a)\rho'(b).$$
This is a statement in the tangent space at $\rho_0(ab)$ in $G$. In order to translate all of the tangent vectors back to the identity, we multiply this equation by $\rho_0(ab)^{-1}$ :
\begin{eqnarray*}\rho'(ab)\rho_0(ab)^{-1} &=& \rho'(a)\rho_0(a)^{-1} + \rho_0(a)\rho'(b)\rho_0(b)^{-1}\rho_0(a)^{-1}\\
&=&\rho'(a)\rho_0(a)^{-1} + Ad_{\rho_0(a)} (\rho'(b)\rho_0(b)^{-1}).
\end{eqnarray*}
Letting $\mathfrak g = T_{\operatorname{Id}} G$ denote the Lie algebra of $G$, define $z : \pi_1 M \rightarrow \mathfrak g$ by $z(\gamma) = \rho'(\gamma)\rho_0(\gamma)^{-1}$. Then $z$ satisfies the \emph{cocycle condition}:
\begin{equation} z(ab) = z(a) + Ad_{\rho_0(a)}z(b). \label{cocycle} \end{equation}
The group cocycles $Z^1(\pi_1 M, \mathfrak g_{Ad \rho_0})$ are defined to be all functions $z$ satisfying Equation~\ref{cocycle} for all $a,b \in \pi_1 M$. The cocyle $z \in Z^1(\pi_1 M, \mathfrak g_{Ad \rho_0})$ is called an \emph{infinitesimal deformation} of the representation $\rho_0$. Next, suppose $\rho_t = g_t \rho_0 g_t^{-1}$ is a trivial deformation of $\rho_0$. Differentiating shows that $$\rho'(\gamma)\rho_0(\gamma)^{-1} = g' - Ad_{\rho_0(\gamma)} g'.$$
The co-boundaries $B^1(\pi_1 M, \mathfrak g_{Ad \rho_0})$ are defined to be all group cocycles $z$ such that $ z(\gamma) = u - Ad_{\rho_0(\gamma)} u$ for some $u \in \mathfrak g$. These are the trivial infinitesimal deformations. Now define the cohomology group $H^1(\pi_1 M, \mathfrak g_{Ad \rho_0}) = Z^1(\pi_1 M, \mathfrak g_{Ad \rho_0}) / B^1(\pi_1 M, \mathfrak g_{Ad \rho_0}).$
If $\mathscr R(\pi_1 M; G)$ is a smooth manifold at $\rho_0$, then $H^1(\pi_1 M, \mathfrak g_{Ad \rho_0})$ describes the tangent space at $\rho_0$.
In all cases of interest here, $G$ is an algebraic group and the representation space $\mathscr R(\pi_1 M; G)$ can be given the structure of an algebraic variety.

\subsection{Projective Geometry}
Real projective geometry will provide the framework for the constructions in this paper. It is more flexible than metric geometries, but nonetheless has a lot of useful structure. The real projective space $\RP^n$ is the space of lines in $\RR^{n+1}$. It is an $n$-dimensional manifold, orientable if and only if $n$ is odd. The group $\GL(n+1, \RR)$ acts by diffeomorphisms on $\RP^n$, with kernel given by its center $\{ \lambda I : \lambda \in \RR^*\}$. Thus $\PGL(n+1,\RR)$, defined to be the quotient of $\GL(n+1,\RR)$ by its center, acts faithfully by diffeomorphisms on $\RP^n$. A hyperplane of dimension $k+1$ in $\RR^{n+1}$ descends to a copy of $\RP^k$ inside $\RP^n$, which we call a \emph{$k$-plane}. The \emph{lines} in $\RP^n$ are described by the case $k=1$. They correspond to two-dimensional planes in $\RR^{n+1}$. Note that $k$-planes in $\RP^n$ are taken to other $k$-planes by $\PGL(n+1,\RR)$, so these are well defined geometric objects in projective geometry; they play the role of totally geodesic hyperplanes in a Riemannian model geometry.

A \emph{projective structure} on a manifold $M^n$ is an $(G,X)$ structure for $G = \PGL(n+1,\RR)$, $X = \RP^n$.
Hyperbolic structures, $\AdS$ structures, and half-pipe structures (to be defined in Section~\ref{HP}) are special examples of projective structures.

\subsection{Hyperbolic geometry and cone singularities}
\label{hyp-cone}

Given a smooth manifold $M$, the data of a hyperbolic metric on $M$ (i.e. a Riemannian metric with sectional curvatures equal to $-1$) is equivalent to a $(G,X)$ structure on $M$, where $X = \HH^n$ and $G = \text{Isom}(\HH^n)$ is the group of isometries of $\HH^n$. For basics of of hyperbolic geometry, see \cite{Thurston, Thurston-book, Ratcliffe-book}. We quickly recall the definitions, and then describe the main structures of interest: hyperbolic manifolds with cone singularities.

Let $\RR^{n,1}$ denote $\RR^{n+1}$ equipped with the $(n,1)$ Minkowski form $\eta$:
\begin{equation*}
\eta = \minimatrix{-1}{0}{0}{I_n}.
\end{equation*}
The projective model for hyperbolic space is given by the negative lines with respect to this quadratic form:
$$\HH^n = \{ x : x^T \eta x < 0\}/\RR^*.$$ The group $\PO(n,1) \subset \PGL(n+1, \RR)$ of matrices (up to $\pm I$) that preserve $\eta$ defines the isometry group in this model. So, every hyperbolic structure is also a projective structure and we say that hyperbolic geometry is a \emph{specialization} or \emph{sub-geometry} of projective geometry. Geodesic lines and hyperplanes in $\HH^n$ are given by lines and hyperplanes in $\RP^n$ that intersect $\HH^n$.

We now restrict to dimension three. Let $N$ be a closed oriented three-manifold, with $\Sigma$ a knot in $N$. Let $M = N \setminus \Sigma$. A \emph{hyperbolic cone-manifold structure} on $(N,\Sigma)$ is given by a smooth hyperbolic structure on $M$ such that the geodesic completion is topologically $N$ (see \cite{Hodgson-98}). The holonomy representation $\rho$ for the cone manifold structure refers to the holonomy representation for the smooth structure on $M$.

\begin{figure}[h]
{
\centering

\def\svgwidth{2.0in}
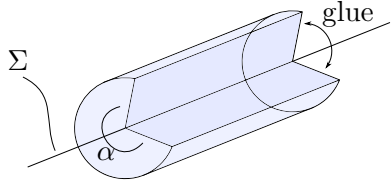

}
\caption[A cone singularity in $\HH^3$.]{Schematic of a cone singularity in $\HH^3$ (or any Riemannian model geometry). Perpendicular to the singular locus, the geometry is that of a cone-point on a surface.\label{cone-H3-glue}} 
\end{figure}

Consider a tubular neighborhood $T$ of $\Sigma$. The developing map $D$ on $\widetilde{T\setminus \Sigma}$ extends to the geodesic completion $\widetilde{T \setminus \Sigma} \cup \widetilde \Sigma$, which is the universal branched cover of $T$ branched over $\Sigma$. The image $D(\widetilde \Sigma)$ is a one-dimensional set in $\HH^3$ which must be fixed point-wise by the holonomy $\rho(\mu)$ of a meridian $\mu$ encircling $\Sigma$. 
We will assume $\rho(\mu)$ is non-trivial. Then $\rho(\mu)$ is a rotation about a geodesic $\Cln$ in $\HH^3$ and $D$ maps $\widetilde \Sigma$ diffeomorphically onto $\Cln$. In this sense, the singular locus $\Sigma$ is totally geodesic.

 The local geometry at points of $\Sigma$ is determined by the \emph{cone angle} $\alpha$, given up to a multiple of $2\pi$, by the rotation angle of $\rho(\mu)$. To determine the full cone angle, one may lift the developing map $D$ to $\widetilde \HH^3 \setminus \Cln$, yeilding a lift of $\rho\big|_{\pi_1 T}$ to the universal cover $\widetilde H_\Cln$ of the subgroup $H_\Cln \subset \PO(n,1)$ that preserves $\Cln$. The lifted holonomy of $\mu$ then measures the full cone angle. Note that totally geodesic discs orthogonal to the singular locus exist and are isometric to $\HH^2$ cones of the same cone angle $\alpha$.

We mention that it is equivalent to describe the geometry of $(N,\Sigma)$ near points of $\Sigma$ using finitely many charts that extend to $\Sigma$. That is, for a point $p \in \Sigma$, there exists simply connected open sets $U_1, \ldots, U_k$ in $M$ and charts $\phi_1, \ldots, \phi_k$, such that the closures $\overline{U}_1, \ldots, \overline{U}_k$ cover a neighborhood of $p$ in $N$ and each chart $\phi_i$ extends continuously to $\overline{U}_i$. 
One may also construct a cone singularity by glueing together one or more polygonal wedges around a line as in Figure~\ref{cone-H3-glue}.

\subsection{$\AdS$ geometry and tachyons}
\label{AdS-tachyon}

Anti de Sitter ($\AdS$) geometry is a Lorentzian analogue of hyperbolic geometry in the sense that $\AdS^n$ has all sectional curvatures equal to $-1$. However, the metric on $\AdS^n$ is indefinite of signature $(n-1,1)$. Vectors of negative length-squared are called \emph{time-like}, vectors of positive length-squared are called \emph{space-like}, and non-zero vectors with zero length are called \emph{light-like} or \emph{null}. For basics on Lorentzian geometry, see \cite{Beem-book} or any physics text on relativity. The implications of negative curvature in Lorentzian geometry are somewhat different than in Riemannian geometry. For example, $\AdS^n$ has an ideal boundary at infinity, but only space-like and light-like geodesics have endpoints on this ideal boundary. Time-like geodesics, on the other hand, are periodic. The geometry in the time-like directions acts more like a positively curved Riemannian space. The reader may consult \cite{Benedetti-09} or \cite{Danciger-11} for some basics of $\AdS$ geometry.

Let $\RR^{n-1,2}$ denote $\RR^{n+1}$ equipped with the $(n-1,2)$ Minkowski form $\eta$, which we choose to write as follows
\begin{equation*}
\eta = \begin{pmatrix} -1 & 0 & 0\\  0 & I_{n-1} & 0\\ 0 & 0 & -1 \end{pmatrix}.
\end{equation*}
The projective model of $\AdS^n$ is given by: 
$$\AdS^n = \{x \in \RR^{n+1} : x^T \eta x < 0 \}/\RR^*.$$
The group $\PO(n-1,2)$ of matrices (up to $\pm I$) that preserve $\eta$ defines the isometry group in this model.
If $n$ is even, then $\PO(n-1,2) \cong \SO(n-1,2)$ has two components, one that preserves time-orientation and one that reverses it.
If $n$ is odd, $\PO(n-1,2)$ has four components corresponding to the binary conditions orientation-preserving (or not), and time-orientation preserving (or not). The orientation preserving, time-orientation preserving subgroup is the component of the identity $\PO_0(n-1,2) = \PSO_0(n-1,2)$. The orientation preserving subgroup is $\PSO(n-1,2)$.
Totally geodesic lines and planes are given by lines and planes in $\RP^n$ which intersect $\AdS^n$. A line is time-like, light-like, or space-like if the signature of the corresponding $2$-plane in $\RR^{n+1}$ is $(0,2)$, degenerate, or $(1,1)$ respectively. Each codimension one space-like plane, meaning positive definite signature, is isometric to a copy of $\HH^{n-1}$.

We now specialize to dimension three and describe the $\AdS$ analogue of a cone singularity: a tachyon. Let $N$ be a closed three-manifold, with $\Sigma$ a knot in $N$. Let $M = N \setminus \Sigma$.
We give the following definition of an $\AdS^3$ manifold with tachyon singularities. Barbot-Bonsante-Schlenker give an equivalent definition in \cite{Barbot-09} as well as a detailed discussion of tachyons and other singularities in $\AdS$.

\begin{Definition}
An \emph{$\AdS^3$ structure on $N$ with a tachyon at $\Sigma$} is given by a smooth $\AdS^3$ structure on $M$ such that the geodesic completion is topologically $N$. The singular locus $\Sigma$ is required to be space-like, and the local future and local past at points of $\Sigma$ must each be connected and non-empty.
\end{Definition}

\begin{figure}[h]
{
\centering

\def\svgwidth{2.5in}
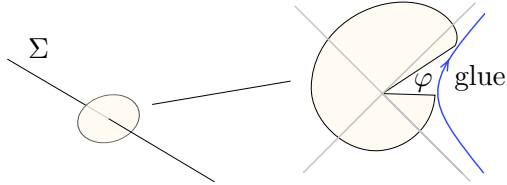

}
\caption[Schematic of a tachyon singularity in $\AdS^3$.]{\label{ugly-fig-AdS} Schematic of a tachyon singularity in $\AdS^3$. The geometry orthogonal to the space-like singular locus can be constructed by glueing a wedge in $\AdS^2$ together with a Lorentz boost. The glueing depicted produces a tachyon of mass $-|\phi|$.} 
\end{figure}

Consider a tubular neighborhood $T$ of $\Sigma$. The developing map $D$ on $\widetilde{T\setminus \Sigma}$ extends to the geodesic completion $\widetilde{T \setminus \Sigma} \cup \widetilde \Sigma$, which is the universal branched cover of $T$ branched over $\Sigma$. The image $D(\widetilde \Sigma)$ is (locally) a one-dimensional set in $\AdS^3$ which must be fixed by the holonomy $\rho(\mu)$ of a meridian $\mu$ around $\Sigma$. 
Assuming $\rho(\mu)$ is non-trivial, it point-wise fixes a geodesic $\Cln$ in $\AdS^3$ and $D$ maps $\widetilde \Sigma$ diffeomorphically onto $\Cln$. In this sense $\Sigma$ is totally geodesic.
By definition, $\Cln$ is required to be \emph{space-like}. The group $G_\Cln$ of orientation preserving isometries that pointwise fix a space-like geodesic $\Cln$ is isomorphic to $\RR\times \ZZ/2$. The $\ZZ/2$ factor is rotation by $\pi$, while the non-compact $\RR$ factor corresponds to a \emph{hyperbolic rotation} around $\Cln$, acting as a Lorentz boost in each tangent plane orthogonal to $\Cln$.  We choose an orientation of $\Cln$ which then determines an orientation on the tangent planes $\Cln^\perp$. The \emph{hyperbolic angle} $\phi$ of a hyperbolic rotation $A \in G_\Cln$ is determined by $$\cosh \phi = v^T\eta Av$$ where $v\in \Cln^\perp$ is any space-like unit vector. The sign of $\phi$ is determined from the orientation of $\Cln^\perp$ by the convention: $\phi > 0$ if $\{v,Av\}$ matches the orientation of $\Cln^\perp$. 
The \emph{tachyon mass} is the hyperbolic angle of $\rho(\mu)$, provided that $\mu$ is chosen to wind around $\Sigma$ in the direction consistent with the chosen orientation of $\Cln$. Note that
the sign of the tachyon mass is independent of the chosen orientation of $\Cln$. The length of $\Sigma$ is given by the translational part of the holonomy $\rho(\lambda)$ of a longitude $\lambda$ of $\Sigma$. In general $\rho(\lambda)$ may also have a component in $G_\Cln$.

As in the Riemannian case, the local geometry at a singularity is, in general, captured by a lifted version of the local holonomy with image in the lift $\widetilde G_\Cln$ of $G_\Cln$ to the universal cover $\widetilde{\AdS^n \setminus \Cln}$. The group $\widetilde G_\Cln \cong \RR \times \pi \ZZ$, is generated (factor-wise) by hyperbolic rotations, and rotations by integer multiples of $\pi$. The $\pi \ZZ$ component of an element of $\widetilde G_\Cln$ is called the \emph{discrete rotational part}.

\begin{Proposition}
The discrete rotational part of the holonomy around $\Sigma$ is $2\pi$.
\end{Proposition}
\begin{proof}
This follows from the condition that the local future and local past at points of $\Sigma$ must be connected and non-empty. Choose a representative $\mu(t)$ for the meridian so that for every $t$, $\mu(t)$ lies on a ray orthogonal to $\Sigma$ emanating from $p$. Then $D(\mu(t))$ lies entirely in the plane $\Cln^\perp$ orthogonal to $\Cln$ at the point $q = D(p)$. As the future of $p$ and the past of $p$ each have one component,
 $D(\mu(t))$ crosses the four light-like rays emanating from $q$ in $\Cln^\perp$ exactly once (counted with sign). 
\end{proof}

\begin{figure}[h]
{
\centering

\def\svgwidth{2.4in}
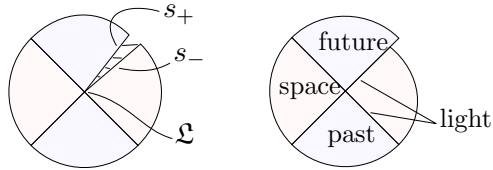

}
\caption[The two-dimensional cross section of a tachyon.]{\label{tachyon-glue-null}A two-dimensional cross-section of a tachyon can be constructed by cutting along a light-like ray and then glueing back together with a Lorentz boost, which act as a dilation along the ray. The figure depicts a tachyon of \emph{negative mass}. This construction should be compared with the construction of Figure~\ref{ugly-fig-AdS}, which produces the same geometry.}
\end{figure}

As with cone singularities, the geometry near a point of the tachyon $\Sigma$ can be equivalently described using finitely many charts that extend to $\Sigma$.
One may also construct tachyons by glueing together one or more polygonal wedges around a space-like line as depicted in Figure~\ref{ugly-fig-AdS}. If the mass is negative this construction can be performed using one wedge with space-like faces as in Figure~\ref{ugly-fig-AdS}.  Interestingly, there is another way to construct a tachyon which is not directly analogous to the Riemannian case. Begin with a tubular neighborhood $U$ of a space-like geodesic $\Cln$ in $\AdS$. The light-cone $\mathcal C$ of $\Cln$ consists of all light-like lines passing through $\Cln$ orthogonally; it is the union of two light-like (degenerate) planes which cross at $\Cln$. Removing $\Cln$ from $\mathcal C$ gives four disjoint open half-planes. Note that the hyperbolic rotations around $\Cln$ preserve each of these light-like half-planes. A tachyon is produced by slitting $U$ along one such half-plane and reglueing via a hyperbolic rotation. Figure~\ref{tachyon-glue-null} depicts a two dimensional cross section of this construction.

\subsection{Transversely hyperbolic foliations}
\label{hyperbolic-foliations}

Let $X$ be a $(n-k)$-dimensional model geometry. A \emph{transversely $(G,X)$ foliation} on a manifold $M^n$ is a smooth foliation by $k$-dimensional leaves so that locally the space of leaves has a $(G,X)$ structure. More concretely, a transversely $(G,X)$ foliation is defined by charts $\phi_\alpha : U_\alpha \rightarrow  \RR^k\times X$ so that each transition map $\phi_\alpha \circ \phi_\beta^{-1} = (f,g)$ respects the product structure and acts on the first factor by a smooth function $f$ and on the second factor by the restriction of an element $g \in G$. As we do not require the smooth functions $f$ to be analytic, a transversely $(G,X)$ foliation is not itself a $(G',X')$ structure.

Consider the case $k=1$, with $X = \HH^{n-1}$, $G = \Isom(\HH^{n-1})$. Then a transversely $(G,X)$ structure on $M$ is called a \emph{transversely hyperbolic foliation}. One can analytically continue charts in the usual way to build a pseudo-developing map $D : \widetilde M \rightarrow X$, which is a local submersion equivariant with respect to a representation $\rho : \pi_1 M \rightarrow G$, again called the holonomy representation. This degenerate developing map encapsulates all of the information about the foliation and its transverse structure. 

Transversely hyperbolic foliations arise as limits of degenerating hyperbolic structures. Assume for simplicity that $M$ is orientable. Consider a path $D_t: \widetilde M \rightarrow \HH^n$ of developing maps for hyperbolic structures such that $D_0 = \lim_{t \rightarrow 0} D_t$ collapses to a local submersion onto a codimension one hyperbolic space $\plane \cong \HH^{n-1}$. The limit $D_0$ will be equivariant with respect to the limiting holonomy representation $\rho_0$. The image of $\rho_0$ must lie in the subgroup $H$ of $\Isom^+(\HH^n)$ that preserves the plane $\plane$. This group $H$ is exactly the isometries of $\plane$, so $D_0$ defines a transversely $(\plane, H) \cong (\HH^{n-1}, \Isom(\HH^{n-1}))$ structure on $M$.

The following theorem of Thurston classifies the topology of closed three-manifolds $M$ that admit a transversely hyperbolic foliation:
\begin{Thm}[Thm 4.9 \cite{Thurston}]
\label{closed-foliation}
Suppose $M^3$ is a closed manifold endowed with a transversely hyperbolic foliation. Let $D$ be a pseudo-developing map with holonomy $\rho$. Then one of the following holds.
\begin{itemize}
\item[(a)] The holonomy group $\rho(\pi_1 M)$ is discrete and $D$ descends to a Seifert fibration $$D_{/\pi_1 M} : M \rightarrow \HH^2/ \rho(\pi_1 M).$$
\item[(b)] The holonomy group $\rho(\pi_1 M)$ is not discrete, and $M$ fibers over the circle with fiber a torus.
\end{itemize}
\end{Thm}

Theorem~\ref{thm:regen-2mm} of the introduction produces hyperbolic cone manifolds that collapse, as the cone angle approaches $2\pi$, to the transversely hyperbolic foliation defined by a Seifert fibration as in case (a) above. 

%
%

\section{Transition theory: half-pipe structures}

\label{HP}
Our description of the transition between hyperbolic and $\AdS$ geometry hinges on the understanding of an interesting new transitional geometry, which we call \emph{half-pipe} or $\HP$ geometry, that bridges the gap between hyperbolic and $\AdS$ geometry. Recall that we wish to construct transitions in the context of hyperbolic and anti de Sitter structures that collapse onto a co-dimension one hyperbolic space. Therefore our model for $\HP^n$ should be the ``midpoint" in a family of models $\XX_s$ which share a common embedded co-dimension one hyperbolic space. We give a natural construction of such a family of models inside of real projective geometry. Though the main focus will be the case $n=3$, we develop this part of the theory in all dimensions $n \geq 2$.

\subsection{$\mathbb H^n$ and $\AdS^n$ as domains in $\mathbb{RP}^n$}
\label{projective-models}

Consider the family $\eta_s$ of diagonal forms on $\mathbb R^{n+1}$ given by $$\eta_s = \begin{pmatrix} -1 & \ 0 \ & 0\\ 0 & I_{n-1} & 0\\ 0 & 0 & \operatorname{sign}(s) s^2 \end{pmatrix}$$ where $s$ is a real parameter and $I_{n-1}$ represents the identity matrix. Each form $\eta_s$ defines a region $\XX_s$ in $\mathbb{RP}^{n}$ by the inequality 
$$x^T \eta_s x = -x_1^2 + x_2^2 + \ldots +x_n^2 + \operatorname{sign}(s) s^2 x_{n+1}^2 < 0.$$
For each $s$, $\XX_s$ is a homogeneous sub-space of $\mathbb{RP}^n$ which is preserved by the group $G_s$ of linear transformations that preserve $\eta_s$.
The usual projective model for hyperbolic geometry is given by $\mathbb H^n = \XX_{+1}$, with $G_{+1} = \PO(n,1)$. In fact, for all $s > 0$ an isomorphism $\XX_{+1} \rightarrow \XX_{s}$ is given by the \emph{rescaling map}
$$\resc_{s} = \minimatrix{I_n}{0}{0}{|s|^{-1}} \ \in \ \PGL(n+1,\mathbb R).$$ Note that $\resc_s$ conjugates $\PO(n,1)$ to $G_s$. Similarly, $\XX_{-1}$ is the usual projective model for anti de Sitter geometry, $\AdS^n$, with $G_{-1} = \PO(n-1,2)$. For all $s < 0$, the map $\resc_s$ gives an isomorphism $\XX_{-1} \rightarrow \XX_{s}$, conjugating $\PO(n-1,2)$ into $G_{s}$. The rescaling map $\resc_s$ should be thought of as a projective change of coordinates which does not change intrinsic geometric properties.
\begin{Remark}
\label{hyperboloid}
For $s \neq 0$, a constant curvature $-1$ metric on $\XX_s$ is obtained by considering the hyperboloid model, defined by $x^T \eta_s x = -1$. In this sense, the maps $\resc_s$ are isometries.
\end{Remark}
\noindent There is a distinguished codimension one totally geodesic hyperbolic space $\plane^{n-1}$ defined by $$x_{n+1}= 0 \ \ \text{ and } -x_1^2 + x_2^2 + \ldots + x_{n}^2 <  0.$$ Note that $\plane^{n-1}$ is contained in $\XX_s$ for all $s$. The rescaling map $\resc_s$ point-wise fixes $\plane^{n-1}$.

\subsection{Rescaling the degeneration - definition of $\HP^n$}
\label{rescaling-defnHP}
\begin{figure}[h]
\centering
\includegraphics[height = 0.8in]{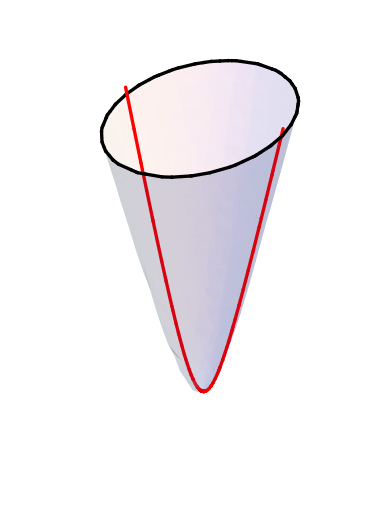}
\includegraphics[height = 0.8in]{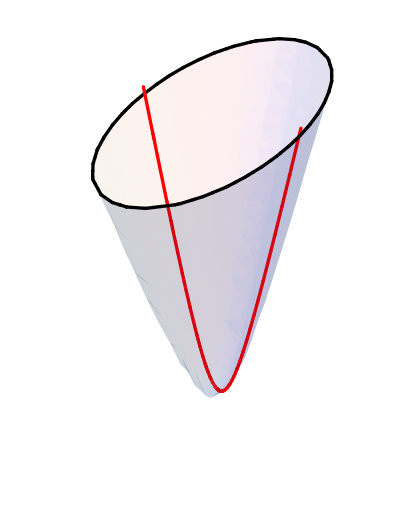}
\includegraphics[height = 0.8in]{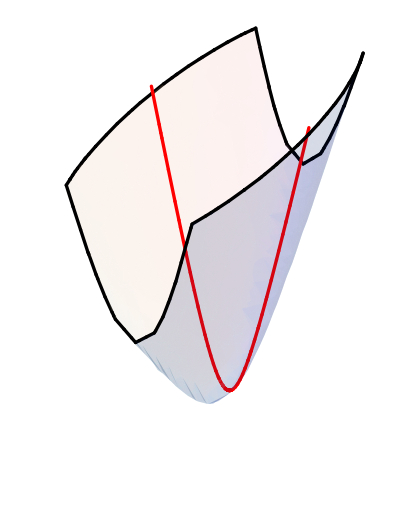}
\includegraphics[height = 0.8in]{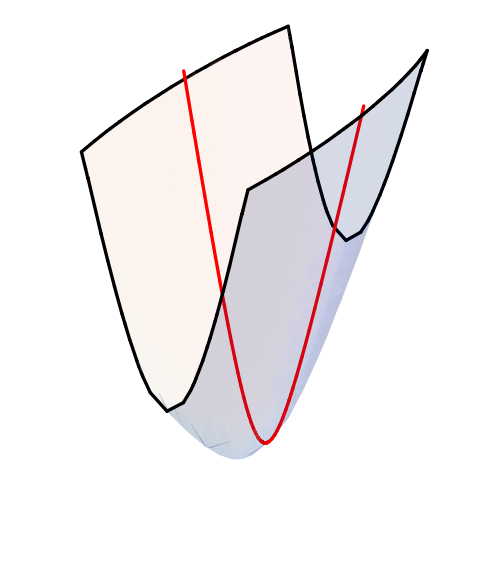}
\includegraphics[height = 0.8in]{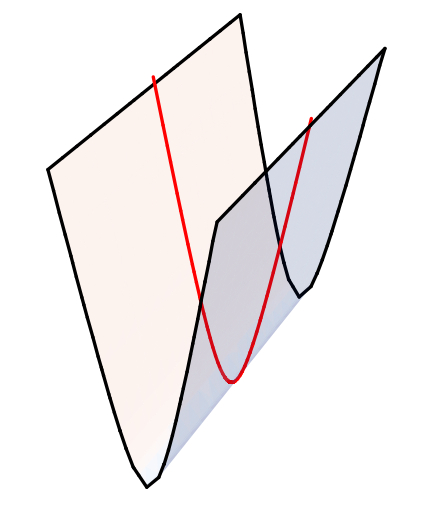}
\caption[Hyperboloid models converging to a degenerate hyperboloid.]{For each $s > 0$, the hyperboloid $x^T\eta_s x = -1$ gives a model for $\HH^2$ (left four figures). As $s \rightarrow 0^+$, the limit is the hyperboloid model for $\HP^2$ (shown right). The distinguished codimension one hyperbolic space $\plane \cong \mathbb H^1$ is shown in red.}
\label{limit-HP}
\end{figure}

The space $\XX_0$ is a natural intermediary space between $\mathbb H^n$ and $\AdS^n$. However, as the metric $\eta_0$ is degenerate, the full group of isometries of $\XX_0$ makes the structure too flimsy to be of much use in our transition context. In order to determine a useful structure group for $\XX_0$ we examine the degeneration context in which we hope to construct a transition. In this section, we will not pay close attention to technical details about collapsing.

Consider a family of developing maps $D_t : \widetilde M \rightarrow \XX_{1}$ with holonomy $\rho_t : \pi_1 M \rightarrow G_{1} = \text{PO}(n,1)$, defined for $t > 0$. Suppose that at time $t=0$, our developing maps collapse to $D_0$, a local submersion onto the co-dimension one hyperbolic space $\plane^{n-1}$. In particular the last coordinate $x_{n+1}$ converges to the zero function. The holonomy representations $\rho_t$ then converge to a representation $\rho_0$ with image in the subgroup $H_0\cong \PO(n-1,1)$ that preserves $\plane^{n-1}$. The one dimensional foliation defined by the local submersion $D_0$ has a natural transverse $\HH^{n-1}$ structure. The foliation together with its transverse structure is called a \emph{transversely hyperbolic foliation} (see Section~\ref{hyperbolic-foliations}). We assume for simplicity that the the fibers of the foliation can be consistently oriented so that in particular the holonomy representation $\rho_0$ of the transverse structure has image in the subgroup 
$$H_0^+ =\left\{ \begin{pmatrix} A & 0\\ 0 & 1 \end{pmatrix} : A \in \SO(n-1,1) \right\}/\{\pm I\} \ \  \cong \ \ \PSO(n-1,1).$$ Next, apply the rescaling map $\resc_{t}$ to get the family $\resc_{t} D_t : \widetilde M \rightarrow \XX_{t}$. This does not change the intrinsic hyperbolic geometry, but extrinsically in $\mathbb{RP}^n$ this stretches out the collapsing direction: $\resc_t$ rescales the $x_{n+1}$ coordinate by $1/t$. Let us assume that $\resc_t D_t$  converges as $t \rightarrow 0$ to a local diffeomorphism $\mathcal D : \widetilde M \rightarrow \XX_0$. The map $\mathcal D$ will be equivariant with respect to a representation $\rho_{\mathcal D}: \pi_1 M \rightarrow \PGL(n+1,\mathbb R)$. This representation is the limit of the holonomy representations for the $\XX_t$ structures determined by $\resc_t D_t$, which are given by the representations $\resc_t \rho_t \resc_t^{-1}$. For a particular $\gamma \in \pi_1 M$, we write $\rho_t(\gamma) = \minimatrix{A(t)}{w(t)}{v(t)}{a(t)}$ where $A$ is $n\times n$, $w, v^T \in \mathbb R^n$, and $a \in \mathbb R$. Then 
\begin{equation} \resc_t\rho_t(\gamma)\resc_t^{-1} = \minimatrix{A(t)}{t w(t)}{\frac{v(t)}{t}}{a(t)} \ \xrightarrow[\ t \rightarrow 0 \ ]{} \ \minimatrix{A(0)}{0}{v'(0)}{1} \ \ = \ \rho_{\mathcal D}(\gamma).
\label{limit-hol}
\end{equation}
The special form of $\rho_\mathcal{D}$ motivates the following definition.
\begin{Definition}
Define $\HP^n = \XX_0$ and $G_{\HP}$ to be the subgroup of $\PGL(n+1, \mathbb R)$ of matrices (defined up to $\pm I$) with the form $\begin{pmatrix} A & 0\\ v & \pm1\end{pmatrix}$ where $A \in \text{O}(n-1,1)$ and $v^T \in \mathbb R^n$. 
We refer to $G_{\HP}$ as the group of \emph{half-pipe isometries}. A structure modeled on $(\HP^n, G_{\HP})$ is called a \emph{half-pipe structure}. 
\end{Definition}

\begin{Definition}
\label{compatibility}
We say that any path of $\text{O}(n,1)$ representations $\rho_t$ satisfying the limit~(\ref{limit-hol}) is \emph{compatible to first order} at $t=0$ with $\rho_\mathcal{D}$. 
\end{Definition}

Both Lie algebras $\mathfrak{so}(n,1)$ and $\mathfrak{so}(n-1,2)$ split with respect to the adjoint action of $\OO(n-1,1)$ as the direct sum $\mathfrak{so}(n-1,1)\oplus\mathbb R^{n-1,1}$. In both cases, the $\mathbb R^{n-1,1}$ factor describes the tangent directions normal to $\OO(n-1,1)$. The group $G_{\HP}$ is really a semi-direct product $$G_{\HP} \cong  \mathbb R^{n-1,1} \rtimes \OO(n-1,1)$$ where an element $\minimatrix{A}{0}{v'}{\pm1}$ is thought of as an infinitesimal deformation $v'$ of the element $A$ normal to $\OO(n-1,1)$ (into either $\OO(n,1)$ or $\OO(n-1,2)$).

We also note that the isotropy group of a point in $\HP^n$ is $\RR^{n-1} \rtimes (\OO(n-1)\times \ZZ/2).$ The subgroup that fixes a point and also preserves orientation and the orientation of the degenerate direction is $\RR^{n-1} \rtimes \SO(n-1)$.

\subsection{Example: singular torus}
\begin{figure}[h]
\centering

\def\svgwidth{4.3in}
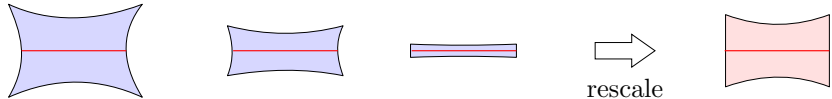

\caption[Collapsing cone tori are rescaled to converge to an $\HP$ torus.]{Fundamental domains for hyperbolic cone tori collapsing to a circle (shown in red). The collapsing structures are rescaled to converge to an $\HP$ structure (right).}
\label{fig:hyptorus}
\end{figure}

We give an illustrative example in dimension $n=2$ of transitioning singular structures on a torus. 
Let $F_2 = \langle a,b\rangle$ be the free group on two generators. For $t > 0$ define the following representations into $G_{+1} = \PO(2,1)$: 
\begin{eqnarray*}\rho_t(a) = \left(\begin{smallmatrix} 3 & 2\sqrt{2} & 0 \\ 2\sqrt{2} & 3 & 0\\ 0 & 0 & 1\end{smallmatrix}\right), & & \rho_t(b) = \left(\begin{smallmatrix} \sqrt{1+t^2} & 0 & t \\ 0 & 1 & 0\\ t & 0 & \sqrt{1+t^2}\end{smallmatrix}\right).
\end{eqnarray*}
For small $t$, the commutator $\rho_t[a,b]$ is elliptic, rotating by an amount $\theta(t) = 2\pi - 2t + O(t^2).$ These representations describe a family of hyperbolic cone tori with cone angle $\theta(t)$. As $t\rightarrow 0$ these tori collapse onto a circle (the geodesic representing $a$).
Next, we rescale this family to produce a limiting half-pipe representation:
\begin{eqnarray*}
\resc_t \rho_t(a)\resc_t^{-1} &=& \rho_t(a) \ \ \ \ \ \ \ \ \ \ \ \ \text{(independent of $t$)}\\
\resc_t \rho_t(b) \resc_t^{-1} &=& \left(\begin{smallmatrix} \sqrt{1+t^2} & 0 & t^2 \\ 0 & 1 & 0\\ 1 & 0 & \sqrt{1+t^2}\end{smallmatrix}\right) \xrightarrow[\ t\rightarrow 0\ ]{} \left(\begin{smallmatrix} 1 & 0 & 0\\ 0 & 1 & 0\\ 1 & 0 & 1\end{smallmatrix}\right).
\end{eqnarray*}
After applying $\resc_t$, the fundamental domains for the hyperbolic cone tori limit to a fundamental domain for a singular $\HP$ structure on the torus (see figure~\ref{fig:hyptorus}). The commutator \begin{equation*}\resc_t \rho_t([a,b])\resc_t^{-1} \xrightarrow[\ t\rightarrow 0\ ]{} \left(\begin{smallmatrix} 1 & 0 & 0\\ 0 & 1 & 0\\ 2 & -2\sqrt{2} & 1\end{smallmatrix}\right)\end{equation*} fixes the singular point and shears in the degenerate direction. This half-pipe isometry can be thought of as an infinitesimal rotation in $\mathbb H^2$.

Next consider the family of singular $\AdS^2$ structures on the torus given by the following $G_{-1} = \PSO(1,2)$ representations defined for $t < 0$: 
\begin{eqnarray*}\sigma_t(a) = \left(\begin{smallmatrix} 3 & 2\sqrt{2} & 0 \\ 2\sqrt{2} & 3 & 0\\ 0 & 0 & 1\end{smallmatrix}\right), & & \sigma_t(b) = \left(\begin{smallmatrix} \sqrt{1-t^2} & 0 & -t \\ 0 & 1 & 0\\ t & 0 & \sqrt{1-t^2}\end{smallmatrix}\right).\end{eqnarray*} Here the commutator $\sigma_t [a, b]$ acts as a Lorentz boost by hyperbolic angle $\phi(t) =  -2t + O(t^2)$ about a fixed point in $\AdS^2$.
These representations describe a family of $\AdS$ tori with a singular point of \emph{hyperbolic angle} $\phi(t)$. The singular point is the Lorentzian analogue of a cone point in Riemannian gometry. We describe the three-dimensional version of this singularity in more detail in Section~\ref{AdS-tachyon}. Again, as $t\rightarrow 0$ these tori collapse onto a circle (the geodesic representing $a$). Similar to the above, we have that 
$\resc_t \sigma_t(b)\resc_t^{-1} \xrightarrow[\ t\rightarrow 0\ ]{} \left(\begin{smallmatrix} 1 & 0 & 0\\ 0 & 1 & 0\\ 1 & 0 & 1\end{smallmatrix}\right)$,
showing that the limiting $\HP$ representation for these collapsing $\AdS$ structures is the same as for the above hyperbolic structures. So we have described a transition on the level of representations. Indeed, applying $\resc_t$ to fundamental domains for the collapsing $\AdS$ structures gives the same limiting $\HP$ structure as in the hyperbolic case above.

\begin{figure}[h]
\centering

\def\svgwidth{4.3in}
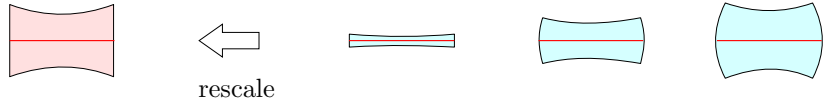

\caption[Collapsing singular $\AdS$ tori limit to the same $\HP$ structure.]{The $\HP$ structure (left) from Figure~\ref{fig:hyptorus} is also the rescaled limit of $\AdS$ tori with ``boost" singularities. Fundamental domains for the $\AdS$ structures are shown right.}
\label{fig:AdStorus}
\end{figure}

\subsection{The geometry of $\HP^{\lowercase{n}}$}
\label{HP-geometry}

Though $\HP^n$ does not have an invariant Riemannian metric, there are some useful geometric tools for studying $\HP$ structures. To begin with, the form $\eta_0$ induces a degenerate metric on $\HP^n$. The degenerate direction defines a foliation of $\HP^n$ by degenerate lines. These are exactly the lines of constant $x_1,\ldots,x_n$ coordinates, with $x_{n+1}$ allowed to vary. There is a projection map $\pi : \HP^n \rightarrow \plane^{n-1} \cong \HH^{n-1}$, given in coordinates by $$\pi(x_1,\ldots,x_n, x_{n+1}) = (x_1,\ldots,x_n,0)$$ which makes the foliation by degenerate lines a (trivial) $\RR$-bundle over $\HH^{n-1}$. The projection commutes with the action of $G_{\HP}$ in the sense that if $g \in G_{\HP}$, then $\pi \circ g = \pi_*(g) \circ \pi$,
where $\pi_*: G_{\HP} \rightarrow \OO(n-1,1)$ is given by $\pi_* \minimatrix{A}{0}{v}{\pm1} = A$. Thus $\pi$ defines a transverse hyperbolic structure on the degenerate lines of $\HP^n$. This transverse structure descends to any $\HP^n$ structure on a manifold $M$. So an $\HP$ structure on $M$ induces a \emph{transversely hyperbolic foliation} on $M$ (see Section~\ref{hyperbolic-foliations}). This can be described directly with developing maps: If $D: \widetilde M \rightarrow \HP^n$ is a local diffeomorphism, equivariant with respect to $\rho: \pi_1 M \rightarrow G_{\HP}$, then $D_0 = \pi \circ D$ is a local submersion onto $\HH^{n-1}$ which is equivariant with respect to $\pi_* \circ \rho : \pi_1 M \rightarrow \OO(n-1,1)$.
Thinking of the induced transversely hyperbolic foliation, we will sometimes refer to the degenerate direction as the \emph{fiber direction}.

Topologically, $\HP^n$ is just $\HH^{n-1} \times \RR$. A particularly useful diffeomorphism is given by $(\pi, L): \HP^n \rightarrow \HH^{n-1} \times \RR$, where $\pi$ is the projection defined above and $L$ is defined in coordinates by $$L(x_1,\ldots,x_n,x_{n+1}) = \frac{x_{n+1}}{x_1\sqrt{1 - (\frac{x_2}{x_1})^2 - \cdots -(\frac{x_n}{x_1})^2}}.$$

\begin{figure}[h]
{\centering

\def\svgwidth{2.3in}
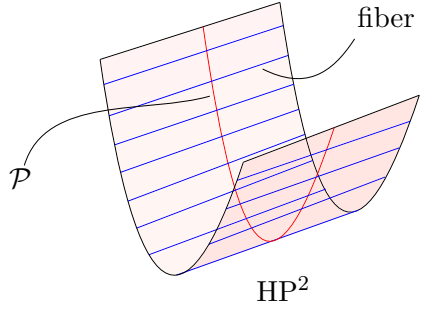

}
\caption[The hyperboloid model of half-pipe geometry.]{The hyperboloid model of half-pipe geometry in dimension two. The degenerate fibers (blue) foliate $\HP^2$.}
\end{figure}
Our choice of structure group $G_{\HP}$ makes the geometry more stiff than the geometry of the degenerate metric alone.
In particular, the non-zero vector field \begin{equation*}
X_{fiber} = x_1\sqrt{1- \left(\frac{x_2}{x_1}\right)^2 - \cdots -\left(\frac{x_n}{x_1}\right)^2} \frac{\partial}{\partial x_{n+1}}
\end{equation*}
descends to $\HP^n \subset \RP^{n}$ and is invariant under $G_{\HP}$ up to $\pm$. It is tangent to the degenerate direction. The group $G_{\HP}$ has four components, corresponding to the conditions orientation preserving (or not) and preserving $X_{fiber}$ (or flipping it). We denote the two components that preserve $X_{fiber}$ by $G_{\HP}^f$. Declaring $X_{fiber}$ to have length one, we can measure lengths along degenerate fibers, as follows. Let $\gamma(t)$ be a path parallel to the degenerate direction, defined for $t \in [a,b]$. Then $\gamma'(t) = f(t)X_{fiber}$
and we define
\begin{eqnarray*}
L_{fiber}(\gamma) &=& \int_a^b f(t)dt.
\end{eqnarray*}
We note that $L_{fiber}(\gamma) = L(\gamma(b)) - L(\gamma(a))$.
\begin{Proposition}
$|L_{fiber}(\gamma)|$ is invariant under $G_{\HP}$. The sign of $L_{fiber}(\gamma)$ is preserved by the subgroup $G_{\HP}^f$ that preserves the fiber direction.
\end{Proposition}
Note that we can not measure such a fiber length for a path transverse to the fiber direction because there is no invariant projection onto the fiber direction. This is the reason that no invariant Riemannian metric exists on $\HP^n$.

\begin{Remark}
Although we have not yet given a detailed discussion of singularities in $\HP$ geometry (see Section~\ref{HP-inf-cone}), we note here that, at least in dimensions $n=2$ and $n=3$, all $\HP$ structures on closed manifolds must have singularities (see the author's thesis~\cite{Danciger-11}) for a proof). This might not be surprising as $\HP$ geometry was designed for the purpose of transitioning from singular hyperbolic structures to singular $\AdS$ structures. 
\end{Remark}

\subsection{Regeneration}
\label{general-regen}
In this section, we show how to regenerate $\mathbb H^n$ and $\AdS^n$ structures from $\HP^n$ structures.

\begin{Proposition}[Regeneration]
Let $M_0$ be a compact $n$-manifold with boundary and let $M$ be a thickening of $M_0$, so that $M\setminus M_0$ is a collar neighborhood of $\partial M_0$. 
Suppose $M$ has an $\HP$ structure defined by developing map $D_{\HP}$, and holonomy representation $\sigma_{\HP}$. 
Let $X$ be either $\mathbb H^n$ or $\AdS^n$ and let $\rho_t: \pi_1 M_0 \rightarrow \Isom(X)$ be a family of representations compatible to first order at time $t = 0$ with $\sigma_{\HP}$ (in the sense of Equation~\ref{limit-hol}). Then we can construct a family of $X$ structures on $M_0$ with holonomy $\rho_t$ for short time.
\label{regen}
\end{Proposition}
\begin{proof}
If $X = \mathbb H^3$, we take $\rho_t$ to be defined for $t \geq 0$, while if $X = \AdS^3$ then we take $\rho_t$ to be defined for $t \leq 0$. This allows us to use the notation from Section~\ref{projective-models} and treat both cases at once.

The representations $\sigma_t := \resc_t \rho_t \resc_t^{-1} : \pi_1 M_0 \rightarrow G_t \subset \PGL(n+1, \mathbb R)$ converge, by assumption, to $\sigma_{\HP}$ in $\PGL(n+1, \mathbb R)$. Thus, thinking of the $\HP$ structure as a projective structure, the Ehresmann-Thurston principle (Proposition~\ref{prop:deform-withboundary}) implies that for small time there is a family of nearby projective structures on $M_0$ with holonomies $\sigma_t$. These projective structures are given by developing maps $F_t : \widetilde M_0 \rightarrow \mathbb{RP}^n$ which converge (in the compact open topology) to $D$ as $t \rightarrow 0$. We show now that $F_t$ is the developing map for an $(G_t, \XX_t)$ structure. We will use the following easy lemma.
\begin{Lemma}
Let $K$ be a compact set  and let $F_t : K \rightarrow \RP^3$ be any continuous family of functions. Suppose $F_0(K)$ is contained in $\XX_s$. Then there is an $\epsilon > 0$ such that $|t| < \epsilon$ and $|r-s| < \epsilon$ implies that $F_t(K)$ is contained in $\XX_r$.
\end{Lemma}
Consider a compact fundamental domain $K \subset \widetilde M_0$. $D(K)$ is a compact set contained in $\HP^n = \XX_0 \subset \mathbb RP^n$. By the lemma, $F_t(K)$ is also contained in $\XX_t$ for all $t$ sufficiently small. Now, since $F_t$ is equivariant with respect to $\sigma_t: \pi_1 M \rightarrow G_t$, we have that the entire image of $F_t$ is contained in $\XX_t$. Thus (for small $t$), $F_t$ determines an $\XX_t$ structure with holonomy $\sigma_t$. Now, applying the inverse of the rescaling map gives developing maps $D_t = \resc_t^{-1} F_t$ into $X$ which are equivariant with respect to $\rho_t = \resc_t^{-1} \sigma_t \resc_t$. These define the desired $X$ structures.
\end{proof}
Note that while this proposition applies in broader generality than Theorem~\ref{thm:main-regen} from the Introduction, it does not guarantee any control of the geometry at the boundary. We study behavior near the boundary in Section~\ref{singular}.

In light of the constructions of this section, we make the following definition of geometric transition. Note that there is no mention of half-pipe geometry in the definition:
\begin{Definition}
Let $M$ be an $n$-dimensional manifold. A \emph{geometric transition} from $\HH^n$ structures to $\AdS^n$ structures is a continuous path of projective structures $\mathscr P_t$ on $M$ so that
\begin{itemize} \item for $t > 0$, $\mathscr P_t$ is conjugate to a hyperbolic structure
\item for $t < 0$, $\mathscr P_t$ is conjugate to an $\AdS$ structure.
\end{itemize}
\end{Definition}

Proposition \ref{regen} implies:
\begin{Proposition}
Let $M$ be a compact manifold with boundary and let $h_t$ be a path of hyperbolic (resp anti de Sitter) structures on $M$ that degenerate to a transverse hyperbolic foliation. Suppose the $h_t$ limit, as projective structures, to an $\HP$ structure. Then a transition to anti de Sitter (resp. hyperbolic) structures can be constructed if and only if the transition can be constructed on the level of representations.
\end{Proposition}

%
%

\section{Singular three dimensional structures}

\label{singular}

In this section, our goal is to build transitions from hyperbolic cone structures to their $\AdS$ analogues, tachyon structures. To do this, we generalize the notion of cone singularity to projective structures.
We focus on dimension three, though much of what is said here applies in general.

\subsection{Cone-like singularities for $\RP^3$ structures}
\label{s:cone-like-3d}

Let $N$ be an orientable three-manifold with $\Sigma \subset N$ a knot. Let $M = N \setminus \Sigma$.

\begin{Definition}
\label{cone-like-3d}
 A \emph{cone-like projective structure with a cone-like singularity} on $(N,\Sigma)$ is a smooth projective structure on $M$ defined by charts $(U_\alpha, \phi_\alpha)$ such that
\begin{itemize}
\item Every chart $\phi_\alpha : U_\alpha \rightarrow \RP^3$ extends continuously to the closure $\overline{U_\alpha}$. In the case $\overline{U_\alpha}$ contains points of $\Sigma$, we require that $\phi_\alpha$ maps $\overline{U_\alpha} \cap \Sigma$ diffeomorphically to a segment of a line $\Cln_\alpha$ in $\RP^3$.
\item For every point $p \in \Sigma$, there is a neighborhood $B$ of $p$ and finitely many charts $(\phi_1,U_1), \ldots, (\phi_k,U_k)$ such that $B$ is covered by $\overline{U_1},\ldots,\overline{U_k}$ and for each $j$, $B \cap \Sigma \subset \overline{U_j} \cap \Sigma.$
\end{itemize}
$M$ is called the \emph{smooth part} and $\Sigma$ is called the \emph{singular locus}. Note that in the case $\overline{U_\alpha \cap U_\beta}$ contains points of $\Sigma$, the transition function $g_{\alpha \beta} \in \PGL(3,\RR)$ maps $\Cln_\beta$ to $\Cln_\alpha$.
\end{Definition}

We note that a projective structure with cone-like singularities on $(N,\Sigma)$ induces an $\RP^1$ structure on $\Sigma$, which is compatible with the projective structure on $M = N \setminus \Sigma$.

\begin{Definition}
Let $(N, \Sigma)$ and $(N', \Sigma')$ be two projective three-manifolds with cone-like singularities. An \emph{isomorphism} $(N, \Sigma) \cong (N', \Sigma')$ is an isomorphism of projective structures $\Phi : N \setminus \Sigma \rightarrow N' \setminus \Sigma'$ which extends to a diffeomorphism $N \rightarrow N'$. We note that $\Phi \big |_\Sigma$ is an isomorphism of the induced $\RP^1$ structures on $\Sigma$ and $\Sigma'$.
\end{Definition}

\begin{Proposition}
\label{model-3cone}
Let $(N, \Sigma)$ be a projective manifold with a cone-like singularity. Let $B$ be a small neighborhood of a point $p \in \Sigma$, with $\Sigma_B = \Sigma \cap B$. Then:
\begin{itemize}
\item The developing map $D$ on $\widetilde{B \setminus \Sigma_B}$ extends to the universal branched cover $\widetilde B = \widetilde{B \setminus \Sigma_B} \cup \Sigma_B$ of $B$ branched over $\Sigma_B$.
\item $D$ maps $\Sigma_B$ diffeomorphically onto an interval of a line $\Cln$ in $\RP^3$.
\item The holonomy $\rho( \pi_1 (B \setminus \Sigma_B))$ point-wise fixes $\Cln$.
\end{itemize}
In particular, there are ``cylindrical'' coordinates $(r,x,y) \in (0,1) \times \RR / \ZZ \times (0,1)$ around $\Sigma_B$ which lift to coordinates on $\widetilde B$ so that $\lim_{r \rightarrow 0} D(r,x,y) =: f(y)$ is a local submersion to $\Cln$ independent of $x$.
\end{Proposition}

\begin{figure}[h]
{\centering

\def\svgwidth{3.6in}
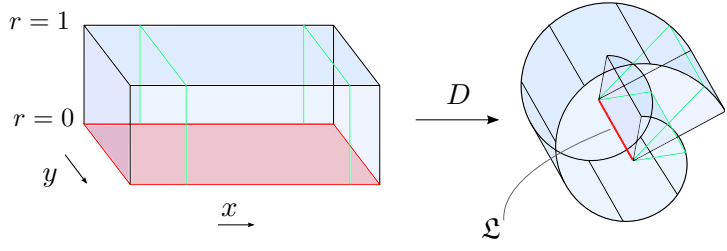

}
\caption{The developing map near a cone-like singularity.}
\end{figure}
\begin{proof}
From the definition of cone-like singularity we may choose $B$ and charts $(\phi_1,U_1),\ldots,$ $(\phi_k, U_k)$ so that $$B = U_1\cup \cdots \cup U_k \cup \Sigma_B$$ and $\cap_{i=1}^k \overline{U_i} = \overline{\Sigma_B}$. By restricting to a smaller neighborhood, we may assume that the following holds:
\begin{itemize}
\item $U_i \cap U_j$ is either empty or $\overline{U_i \cap U_j} \cap \Sigma = \overline{\Sigma_B}$.
\item The $U_i$ are arranged in order around $\Sigma$.
\end{itemize} 
We construct the developing map using $U_1,\ldots, U_k$ by first lifting $U_1$ to $\widetilde{B \setminus \Sigma}$ and mapping into $\RP^3$ with $\phi_1$. Then, the usual analytic continuation process defines $D$ on the rest of $\widetilde{B \setminus \Sigma}$. Note that, by our assumptions $\overline{U_i \cap U_{i+1}} \cap \Sigma$ must be non-empty, and so the transition function $g_{i,i+1}$ maps $\Cln_{i+1}$ to $\Cln_i$. 
Since $\Sigma_B \in \overline{U_i}$ for all $i$, $D$ extends continuously mapping $\Sigma_B$ to $\Cln_1$ by a diffeomorphism. Hence $D$ extends to the universal branched cover $\widetilde B = \widetilde{B \setminus \Sigma_B} \cup \Sigma_B$. The coordinates $(r,x,y)$ are easily obtained by pulling back any choice of cylindrical coordinates around $\Cln_1$ in $\RP^3$.
\end{proof}

The converse of the proposition is also true: Suppose there is a projective structure on $M$  and assume that $\Sigma$ is covered by neighborhoods $B$ so that the developing map $D$ on $\widetilde{B \setminus \Sigma}$ extends to the universal branched cover $\widetilde B$, mapping $B \cap \Sigma$ diffeomorphically to a line. Then it is easy to construct charts satisfying the requirements of Definition~\ref{cone-like-3d}. So $(N,\Sigma)$ is a projective structure with a cone-like singularity. Therefore the two singularities of interest in hyperbolic and $\AdS$ geometry are both examples:

\begin{Proposition}
Let $N$ be a three manifold and $\Sigma$ a knot in $N$. Then
\begin{itemize}
\item \emph{(cone singularities are cone-like):}
The underlying projective structure of a hyperbolic cone structure (see Section~\ref{hyp-cone}) on $(N,\Sigma)$ has a cone-like singularity at $\Sigma$. 
\item \emph{(tachyons are cone-like):}
The underlying projective structure of an $\AdS$ structure with a tachyon (see Section~\ref{AdS-tachyon}) on $(N,\Sigma)$ has a cone-like singularity at $\Sigma$. 
\end{itemize}
\end{Proposition}


Let $G_\Cln$ denote the elements of $\PGL(4,\RR)$ which point-wise fix $\Cln$ and preserve orientation. We fix an orientation of $\Cln$ and an orientation of $\RP^3$ which determines a positive direction of rotation around $\Cln$.  We define the \emph{rotation angle} map $\Rot : G_{\Cln} \rightarrow S^1$ as follows. Given $[A] \in G_{\Cln}$, there is a representative $A$ so that the eigenvalues corresponding to $\Cln$ are both one. Let $\lambda_3, \lambda_4$  be the other eigenvalues. If $\lambda_3, \lambda_4  = \lambda, \bar{\lambda}$ are complex, then $A$ is similar in $\SL(4,\RR)$ to the block diagonal form $$A = \begin{pmatrix} I_2 & 0\\ 0 & |\lambda| R(\theta) \end{pmatrix},$$ where $R(\theta)$ rotates by angle $\theta$ in the positive direction. In this case define $\Rot(A) = e^{i\theta}.$
If $\lambda_3, \lambda_4$ are real, then they both have the same sign and we define $\Rot (A) = \text{sign}(\lambda_2).$ The rotation angle function $\Rot : G_{\Cln} \rightarrow S^1$ is a homotopy equivalence.

Now, consider $D$ and $\rho$ as in Proposition~\ref{model-3cone} and let $\gamma(t) \in \pi_1 (B \setminus \Cln)$. We can find a path $g(t) \in G_\Cln$ such that $D(\gamma(t)) = g(t)\cdot D(\gamma(0))$ with $g(0) = 1$ and $g(1) =\rho(\gamma)$. The path $g(t)$ is unique up to homotopy and defines the \emph{lifted holonomy} $\widetilde \rho(\gamma) \in \widetilde G_\Cln$. Let $m$ be a meridian encircling $\Sigma$ in the direction consistent with the local orientation of $\Sigma$. The geometry in a neighborhood of a point of $\Sigma$ is determined by the lifted holonomy $\widetilde \rho(m)$. There is extra information contained in the lifted holonomy $\widetilde \rho(m)$ that is missing from $\rho(m)$: $\rho(m)$ does not detect how many times $D(m)$ winds around $\Cln$. This information is contained in the \emph{total rotational part} of $\gamma$ defined by the lifted rotation angle map $\widetilde{\mathcal R}: \widetilde G_\Cln \rightarrow \RR$: $$\widetilde{\mathcal R}(\gamma) := \widetilde{\mathcal R} ([g(t)]).$$
The map $\widetilde{\mathcal R} : \pi_1 B \setminus \Sigma \rightarrow \RR$ is a homomorphism. Note that it does not in general extend to a representation of $\pi_1 M$. 
 \begin{Definition}
 \label{rotational-part-3d}
 The quantity $\alpha :=  \widetilde \Rot(m)$ is the \emph{rotational part} of the holonomy at $\Cln$. Note that the rotational part of the holonomy must satisfy $ e^{i\alpha} = \Rot(\rho(\gamma)).$ In the case that the eigenvalues of $\rho(\gamma)$ are real, $\alpha$ is an integer multiple of $\pi$ and we call $\alpha$ the \emph{discrete rotational part} of the holonomy at $\Cln$.
 \end{Definition}


Note that in the case of hyperbolic cone manifolds, the rotational part of the holonomy at $\Sigma$ is exactly the cone angle and determines the local geometry entirely. However, in this more general setting, there can be many geometrically different cones with the same rotational holonomy.

A projective structure with cone-like singularities along a multiple component link $\Sigma$ can be defined analogously. Over the next few sections we will assume $\Sigma$ has one component; this will be the case in the main theorem we are heading towards, and it also keeps the discussion tidy. However, all of the basic theory we develop can easily be extended to the multiple component case.

\subsection{Infinitesimal cone singularities in $\HP^3$}
\label{HP-inf-cone}

In order to develop a theory of geometric transitions with singularities, we consider $\HP$ structures with a singularity that is cone-like with respect to the underlying projective structure. These singularities arise naturally as rescaled limits of collapsing neighborhoods of cone singularities (resp. tachyons) in $\HH^3$ (resp. $\AdS^3$). 

\begin{Definition}
Let $N$ be an oriented three-manifold with $\Sigma \subset N$ a knot. Let $M = N \setminus \Sigma$. An $\HP$ structure with \emph{infinitesimal cone singularity} on $(N, \Sigma)$ is a smooth $\HP$ structure on $M$ whose underlying projective structure has a cone-like singularity at $\Sigma$. Further, we require that there are exactly two degenerate rays emanating from each point of $\Sigma$. Hence $\Sigma$ is a non-degenerate line and the discrete rotational part of the holonomy around $\Sigma$ is $2\pi$.
\end{Definition}

In this section, we describe model neighborhoods around an infinitesimal cone singularity using the $\HP$ geometry. We will show that the local geometry of any infinitesimal cone singularity is realized as the rescaled limit of a model collapsing neighborhood of a cone (resp. tachyon) singularity in hyperbolic (resp. $\AdS$) geometry. We begin by demonstrating this on the level of holonomy representations.

Let $T$ be a solid torus with core curve $\Sigma$ and assume that $T$ has an $\HP$ structure with infinitesimal cone singularity at $\Sigma$. Let $m$ be a meridian encircling $\Sigma$ in the positive direction with respect to the orientation of $\Sigma$. If the holonomy $\rho(m)$ is trivial, then the $\HP$ structure extends smoothly over $\Sigma$, i.e. there is no singularity. This follows from the requirement that the rotational part of the holonomy around $\Sigma$ be $2\pi$. So, we assume that $\rho(m)$ is non-trivial. Then $\rho(m)$ lies in the group $K_\Cln$
of $\HP$ isometries that pointwise fix a non-degenerate line $\Cln$ and preserve both orientation and the direction along degenerate fibers. The holonomy $\rho(\ell)$ of a longitude $\ell$ will lie in the group $H_\Cln$ of $\HP$ isometries that preserve $\Cln$, the orientation of $\Cln$, and the orientation of $\HP^3$.

\begin{figure}[h]
{
\centering

\def\svgwidth{3.0in}
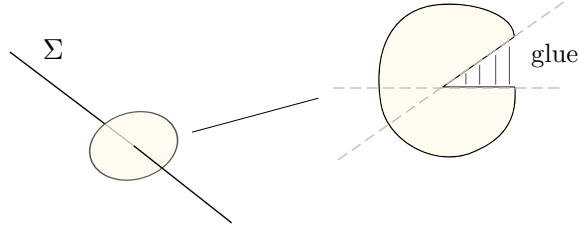

}
\caption[Schematic of an infinitesimal cone singularity in $\HP^3$.]{ \label{inf-cone-wedge} Orthogonal to the singular locus, the geometry can be constructed by glueing together the two non-degenerate boundary rays of a wedge with an infinitesimal rotation.} 
\end{figure}

By conjugating in $G_{\HP}$, we may assume $\rho(m)$ and $\rho(\ell)$ have the following forms: $$\rho(m) = \begin{pmatrix} 1 & 0 & 0 & 0\\ 0 & 1 & 0 & 0\\ 0 & 0 & 1 & 0 \\ 0 & 0 & \omega & 1
\end{pmatrix},
\rho(\ell) =\begin{pmatrix} \cosh d & \sinh d & 0 & 0\\ \sinh d & \cosh d & 0 & 0\\ 0 & 0 & \pm 1 & 0 \\ 0 & 0 & \mu & \pm 1
\end{pmatrix}.$$ The general form of $\rho(m)$ describes $K_\Cln \cong \RR_\omega$, while the general form of $\rho(\ell)$ describes $H_\Cln \cong \RR_d \times \RR_\mu \times \ZZ /2$. The $\RR_d$ factor consists of pure translations along $\Cln$ and the $\ZZ/2$ factor is a rotation by $\pi$ around $\Cln$ which reverses direction along degenerate fibers. We will see how to interpret the $\RR_\mu$ factor shortly. 
Recall that there is a hyperbolic plane $\plane^2 \subset \HP^3$, which we think of as simultaneously lying in each of our family of models $\XX_s$ (refer to Section~\ref{projective-models} for notation). Note that if $\rho(m),\rho(\ell)$ are in the form given above, then the preserved line $\Cln$ lies in $\plane$. If $\rho$ is the limit of rescaled $\PO(3,1)$ representations $\rho_t$, then assuming that $\Cln$ lies in $\plane$ corresponds to assuming that the axis of $\rho_t(m),\rho_t(\ell)$ in $\HH^3$ lies in $\plane$ (at least to first order). Without loss in generality we will assume this throughout the section.

It is easy to construct a path $\rho_t: \langle m, \ell \rangle \rightarrow \PO(3,1)$ whose rescaled limit agrees with $\rho$. Define the path as follows: 
\begin{equation*} \rho_t(m) = \begin{pmatrix} 1 & 0 & 0 & 0\\ 0 & 1 & 0 & 0\\ 0 & 0 & \cos \omega t & -\sin \omega t\\ 0 & 0 & \sin \omega t & \cos \omega t
\end{pmatrix}, \rho_t(\ell) = \begin{pmatrix} \cosh d & \sinh d & 0 & 0\\ \sinh d & \cosh d & 0 & 0\\ 0 & 0 & \pm \cos \mu t & -\sin \mu t\\ 0 & 0 & \sin \mu t & \pm \cos \mu t
\end{pmatrix}.
\end{equation*}
These representations describe hyperbolic cone structures on a tubular neighborhood of $\Sigma$ with cone angles approaching $2\pi$. One easily checks that conjugating $\rho_t$ by the rescaling map $\resc_t$ produces the desired limit as $t \rightarrow 0$. For example:
\begin{equation*}
\resc_t \rho_t(m) \resc_t^{-1} = \begin{pmatrix} 1 & 0 & 0 & 0\\ 0 & 1 & 0 & 0\\ 0 & 0 & \cos \omega t  & -t\sin \omega t \\ 0 & 0 & \sin \omega t/ t & \cos \omega t
\end{pmatrix} \ \xrightarrow[t\rightarrow 0]{} \begin{pmatrix} 1 & 0 & 0 & 0\\ 0 & 1 & 0 & 0\\ 0 & 0 & 1 & 0 \\ 0 & 0 & \omega & 1
\end{pmatrix}.
\end{equation*}

The quantity $\omega$ describes the first order change in rotation angle of $\rho_t(m)$ at $t=0$. Hence we call  $\rho(m)$ an \emph{infinitesimal rotation}. We note that if $\omega > 0$, the cone angle of nearby hyperbolic cone structures must be increasing, while if $\omega < 0$, the cone angle of nearby hyperbolic structures will be decreasing.

\begin{Definition}
The \emph{infinitesimal cone angle} around $\Sigma$ is defined to be the quantity $\omega$. Note that the sign is well-defined and that the lifted holonomy $\widetilde \rho(m)$ is a rotation by $2\pi$ plus an infinitesimal rotation by $\omega$.
\end{Definition}

\begin{Remark}
By an argument using the Schlafli formula, collapsing hyperbolic cone manifolds must have \emph{increasing} cone angle (see \cite{Hodgson-86, Porti-98}).  For this reason, we speculate that singular $\HP$ structures on closed three-manifolds with positive infinitesimal cone angle $\omega > 0$ should not exist in most cases.
\end{Remark}

It is just as easy to construct a path of representations $\rho_t : \langle m, \ell \rangle \rightarrow \PO(2,2)$
whose rescaled limit agrees with $\rho$. Define the path as follows: 
\begin{equation*} \rho_t(m) = \begin{pmatrix} 1 & 0 & 0 & 0\\ 0 & 1 & 0 & 0\\ 0 & 0 & \cosh \omega t & \sinh \omega t\\ 0 & 0 & \sinh \omega t & \cosh \omega t
\end{pmatrix}, \rho_t(\ell) = \begin{pmatrix} \cosh d & \sinh d & 0 & 0\\ \sinh d & \cosh d & 0 & 0\\ 0 & 0 & \pm \cosh \mu t & \sinh \mu t\\ 0 & 0 & \sinh \mu t & \pm \cosh \mu t
\end{pmatrix}.
\end{equation*}
These representations describe $\AdS$ structures on a tubular neighborhood of $\Sigma$ with a tachyon at $\Sigma$ of mass $\omega t$. One easily checks that conjugating $\rho_t$ by the rescaling map $\resc_t$ produces the desired limit as $t \rightarrow 0$. Hence, the infinitesimal angle $\omega$ can also be thought of as an infinitesimal tachyon mass.

Next, we work directly with the $\HP$ geometry at $\Sigma$. Let $p \in \Sigma$ and consider a neighborhood $B$ of $p$. The developing map $D$ on $B \setminus \Sigma$ extends to the universal branched cover $\widetilde B$, branched over $B \cap \Sigma$. The image of $B \cap \Sigma$ is a segment of a non-degenerate line $\Cln$, which we may assume lies in $\plane$. Consider the plane $P$ orthogonal to $\Cln$ and passing through $D(p)$. As $P$ is spanned by the fiber direction and a non-degenerate direction orthogonal to $\Cln$, the restricted metric is degenerate on $P$. The inverse image $C_p := D^{-1}(P)$ is a disk in $B$. Away from $p$, $C_p$ is locally modeled on $\HP^2$. The singularity at $p$ is a cone-like $\HP^2$ singularity. We may, as in the $\HH^3$ and $\AdS^3$ cases, parallel translate $C_p$ (or at least a smaller neighborhood of $p$ in $C_p$) along the interval $\interval = B \cap \Sigma$, giving the identification $B = C_p \times \interval$ near $\interval$. Let $\WW$ be a wedge in $C_p$ (modeled on a wedge in $\HP^2$), and define the product wedge $U = \WW \times \interval$. Product wedges are, as in the hyperbolic and $\AdS$ case, the most natural geometric charts at the singular locus. 

\begin{figure}[h]
{
\centering

\def\svgwidth{2.8in}
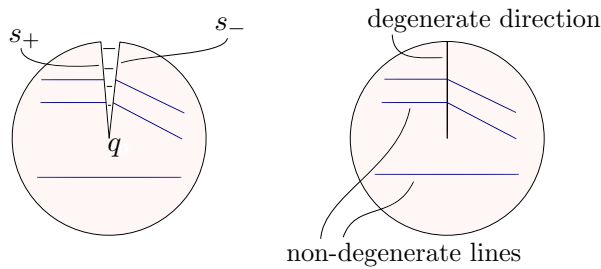

}
\caption[The cross section of an infinitesimal cone singularity.]{ \label{HP-slit-disk-construction} A disk $V$ is slit along a ray in the degenerate direction. It is then glued back together using a non-trivial infinitesimal rotation to produce an \emph{infinitesimal cone singularity}. This construction should be compared with the construction of Figure~\ref{inf-cone-wedge}, which produces the same geometry.}
\end{figure}

We now construct some particularly useful wedges. For simplicity, this part of the discussion will take place in dimension two. The corresponding three-dimensional behavior is easily described by taking the product with a non-degenerate geodesic. Consider the $\HP^2$ cone $C_p$ defined above. By assumption, there are two degenerate rays emanating from $p$.  Pick one of these rays, $r$, and let $\WW$ be $C_p$, but with a slit along the ray $r$, so that the boundary of $\WW$ contains two copies $r_+, r_-$ of $r$ with opposite orientation. Though it is a bit of an abuse, we count $\WW$ as a wedge. It is isomorphic to a disk $V$ in $\HP^2$ with a slit along a degenerate ray $s$ emanating from the center $q$ of $V$. The boundary of $V$ contains two copies $s_+,s_-$ of $s$. We take $s_+$ to be the positive ray, meaning that it is adjacent to the portion of $V$ containing a small positive rotation of $s$. The glueing map $g$ identifies $s_+$ to $s_-$ by an infinitesimal rotation fixing $q$. Note that $g$ fixes $s$ point-wise. Nonetheless, the geometry at $q$ is singular, for the glueing map does not preserve the lines transverse to $s$ (see Figure~\ref{HP-slit-disk-construction}). The holonomy around $p$ is a rotation by $2\pi$ composed with $g$.

\begin{figure}[h]
{
\centering

\def\svgwidth{3.6in}
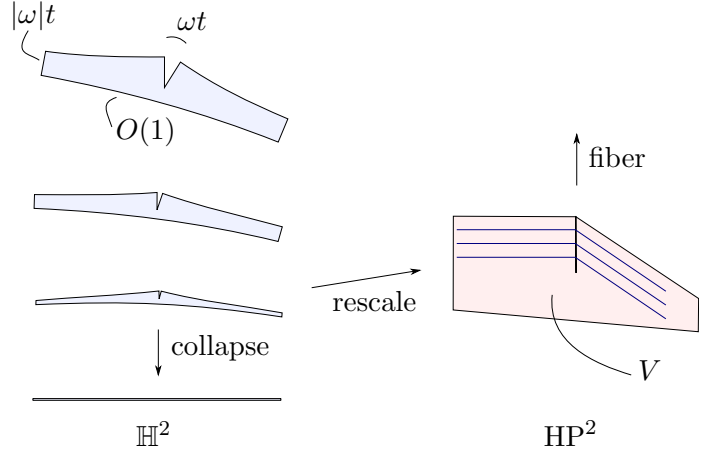

}
\caption[Collapsing hyperbolic cones converge to an $\HP$ cone.]{ \label{model-degen-H2} Polygonal hyperbolic wedges collapse onto a line as the (interior) wedge angle $2\pi - |\omega| t$ approaches $2\pi$. Each wedge is glued together to form a rectangle with a cone point at the center. The rescaled limit of the wedges $V$ is an $\HP$ polygon with a slit along the fiber direction. Glueing the slit together with the rescaled limit of the glueing maps produces an $\HP$ infinitesimal cone singularity with infinitesimal angle $\omega$. }
\end{figure}

Next, we construct a model degeneration of hyperbolic cones (see Section~\ref{hyp-cone}) which when rescaled converge to a given $\HP$ cone. Again, we give the construction in two dimensions; the three-dimensional case is described by taking the product with a geodesic. We will assume that the infinitesimal cone angle $\omega < 0$, so that we can easily draw a picture. Let $\theta(t) = 2\pi - |\omega| t$. Construct a polygonal wedge $V(t)$ in $\HH^2$ with seven sides, six right angles and a seventh (concave) angle $\theta$ at the center point of the wedge as in Figure~\ref{model-degen-H2}. Glueing $V(t)$ together along the sides adjacent to the center point produces a rectangle with a cone point at the center. We arrange for $V(t)$ to be long and skinny, with width roughly equal to one, and thickness $|\omega | t + O(t^2)$. Further, we arrange one ray $s_+$ of the concave part of the wedge to be aligned with the collapsing direction. The glueing map $g(t)$ is a rotation by $\omega t$. Now, the rescaled limit of these collapsing wedges $V(t)$ produces an $\HP$ wedge $V$ of the type described in the previous paragraph. The glueing map $g$ is the rescaled limit of a rotation by $\omega t$, which is an infinitesimal rotation by $\omega$ (as demonstrated explicitly above).

\begin{figure}[h]
{
\centering

\def\svgwidth{4.0in}
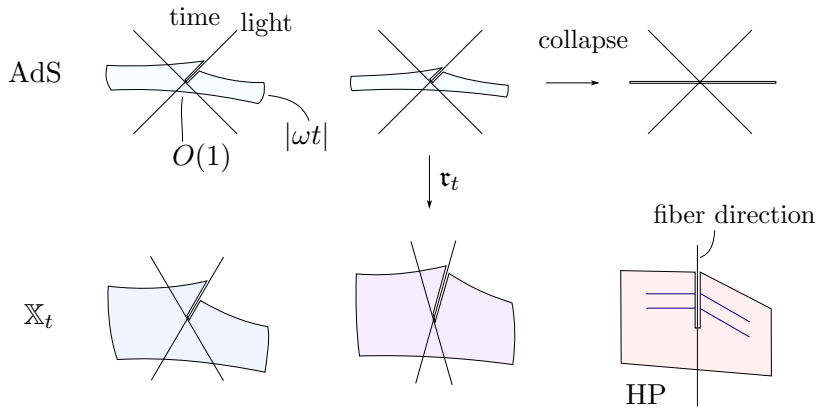

}
\caption[Collapsing $\AdS$ cones collapse to an $\HP$ cone.]{\label{model-degen-AdS}Polygons with a slit along a light-like ray in $\AdS^2$ are glued together with a Lorentz boost of hyperbolic angle $\phi(t) = \omega t$ to form rectangles with a singular point at the center. After rescaling the collapsing time-like direction, these polygons converge to an $\HP$ polygon with a slit along a degenerate ray. This ``wedge" is glued together with the rescaled limit of the Lorentz boosts: an infinitesimal rotation (thought of as an infinitesimal boost) by angle $\omega$. }
\end{figure}

Next, we construct a model degeneration of $\AdS$ tachyons (see Section~\ref{AdS-tachyon}) that when rescaled converge to the given $\HP$ cone. Let $\phi(t) = \omega t$. Let $V(t)$ be a wedge in $\AdS^2$ bounded by seven edges as in Figure~\ref{model-degen-AdS}. The five edges along the convex part of the perimeter should alternate space-like, time-like, space-like, time-like, space-like meeting at four right angles. We arrange for the space-like edges to be of roughly constant length, while the time-like edges have (time-like) length $|\omega | t + O(t^2)$. The two remaining edges $s_+, s_-$ border a slit along a light-like ray emanating from the center $q$ of the wedge. The glueing map $g(t)$, which a is Lorentz boost of hyperbolic angle $\phi$, identifies $s_-$ with $s_+$; the action of $g(t)$ on $s_-$ is a dilation by $e^{\phi}$. Now, the rescaled limit of these collapsing wedges $V(t)$ produces an $\HP$ wedge $V$ of the type described in the previous paragraph. The glueing map $g$ is the rescaled limit of a boost by hyperbolic angle $\omega t$, which is an infinitesimal rotation by $\omega$, alternatively thought of as an infinitesimal boost.

\subsection{Deforming cone-like projective structures}
In order to prove a regeneration theorem for hyperbolic cone (resp. $\AdS$ tachyon) structures, we need to extend the Thurston Ehresmann principle (Proposition~\ref{prop:deform-withboundary}) to the class of projective structures with cone-like singularities. That is we must show that an appropriate deformation of the holonomy representation of a cone-like projective structure produces nearby cone-like projective structures.

Let $N$ be a three-manifold, with $\Sigma \subset N$ a knot, and let $M = N \setminus \Sigma$. Let $T \subset M$ be a neighborhood of $\partial M$ (so $T$ is the result of removing $\Sigma$ from a tubular neighborhood of $\Sigma$ in $N$). Let $\widetilde T$ be the universal cover of $T$. We assume that $\pi_1 T \hookrightarrow \pi_1 M$, so that $\widetilde T$ embeds in $\widetilde M$. The fundamental group $\pi_1 T \cong \mathbb Z \times \mathbb Z$ is generated by the meridian $m$ around $\Sigma$ and a longitude $\ell$.
\begin{Remark}
\label{pi1-injects}
The assumption $\pi_1 T \hookrightarrow \pi_1 M$ holds in every application that we are interested in. However, the assumption is not necessary. Everything done below can be easily modified if a longitude $\ell \mapsto 1$ in $\pi_1 M$.
\end{Remark}

Suppose $(N, \Sigma)$ has a projective structure with cone-like singularity. Let $D: \widetilde T \rightarrow \RP^3$ be the developing map on a chosen lift $\widetilde T$ of $T$, and let $\rho : \pi_1 T \rightarrow \PGL(4,\RR)$ be the holonomy. Using Proposition~\ref{model-3cone} we can construct convenient coordinates $(r,x,y) \in (0,1) \times \RR \times \RR$ for $\widetilde T$ with the following properties:
\begin{itemize}
\item The action of $\pi_1 T$ by deck translations is given by
\begin{align*}
m &: (r,x,y) \mapsto (r,x+1,y) &
\ell : (r,x,y) \mapsto (r,x,y+1).
\end{align*}
\item The limit $\lim_{r \rightarrow 0} D(r,x,y) =: f(y),$ is a local submersion, independent of $x$, to a line $\Cln$ in $\RP^3$. The line $\Cln$ represents the lift of $\Sigma$ corresponding to the chosen lift of $T$.
\item $\rho(m)$ point-wise fixes $\Cln$ and $\rho(\ell)$ preserves (but does not fix) $\Cln$.
\end{itemize}
These coordinates will be useful for proving the following proposition.
\begin{Proposition}
\label{deform-cone-like}
Suppose $\rho_t: \pi_1 M \rightarrow \PGL(4,\RR)$ is a path of representations such that 
\begin{enumerate} \item $\rho_0$ is the holonomy representation of a projective structure with cone-like singularities on $(N,\Sigma)$. Let $\Cln$ be the line in $\RP^3$ fixed by $\rho(\pi_1 \partial M)$. 
\item $\rho_t(m)$ point-wise fixes a line $\Cln_t$, with $\Cln_t \rightarrow \Cln$.
\end{enumerate}
Then, for all $t$ sufficiently small, $\rho_t$ is the holonomy representation for a projective structure with cone-like singularities on $(N,\Sigma)$. 
\end{Proposition}
\begin{proof}
 First, we let $D_0 : \widetilde M \rightarrow \RP^3$ denote the developing map of our projective structure at time $t=0$. Let $M_0 \subset M$ be the result of removing a smaller tubular neighborhood $T' \subset T$ of $\Sigma$ from $M$, so that $M_0$ and $T$ overlap in a neighborhood of $\partial M_0$. By Proposition~\ref{prop:deform-withboundary}, we can deform the projective structure on $M_0$ to get developing maps $D_t : \widetilde {M_0} \rightarrow \RP^3$ that are equivariant with respect to $\rho_t$. Further, by the well-known proof of the proposition, we may assume that $D_t$ converges uniformly in the $\mathscr C^1$ topology on compacts in $\widetilde M_0$. Now we must extend $D_t$ to the rest of $\widetilde M$.

We may assume, by conjugating $\rho_t$ in $\PGL(4,\mathbb R)$, that $\rho_t(m)$ also fixes $\Cln$, for all $t$.
That is, we assume $\Cln_t = \Cln$. 
In order to define $D_t$ on $\widetilde T$ we will need a quick lemma.
\begin{Lemma}
For each $\gamma \in \pi_1 T$, we can take arbitrary powers $\rho_t(\gamma)^z$ in a way that depends smoothly on $z,t$. 
\end{Lemma}
\begin{proof}[Proof of Lemma]
First $\rho_t(\pi_1 T) \subset \PSL(4,\RR)$. So for each $\gamma \in \pi_1 T$, we can find a path $g(t)$ in $\PSL(4,\RR)$ with $g(0) = \text{Id}$ and $g(1) = \rho_0(\gamma)$. The $\log$ function is well-defined sufficiently near to the identity and can be defined by analytic continuation along the path $g(s)$ (this amounts to choosing a branch of $\log$ for the eigenvalues; note that we can not have an odd number of negative real eigenvalues). Next, analytically continue $\log$ along the path $\rho_t(\gamma)$. Thus $\rho_t(\gamma)^z := \exp (z \log \rho_t(\gamma))$ depends smoothly on $z,t$.
\end{proof}
Next, using the coordinates defined above, define $D_t$ on $\widetilde T$ as follows:
\begin{equation*}
D_t(r,x,y) = \rho_t(m)^x \rho_t(\ell)^y \rho_0(m)^{-x} \rho_0(\ell)^{-y} D_0(r,x,y)
\end{equation*}
We check that
\begin{eqnarray*}
D_t(r,x+k,y+j) &=& \rho_t(m)^{x+k} \rho_t(\ell)^{y+j} \rho_0(m)^{-x-k} \rho_0(\ell)^{-y-j} D_0(r,x+k,y+j)\\
&=& \rho_t(m)^{x+k} \rho_t(\ell)^{y+j} \rho_0(m)^{-x-k} \rho_0(\ell)^{-y-j}\rho_0(m)^k\rho_0(\ell)^j D_0(r,x,y)\\
&=& \rho_t(m)^{x+k} \rho_t(\ell)^{y+j} \rho_0(m)^{-x} \rho_0(\ell)^{-y} D_0(r,x,y)\\
&=& \rho_t(m)^{k} \rho_t(\ell)^j \rho_t(m)^{x} \rho_t(\ell)^{y} \rho_0(m)^{-x} \rho_0(\ell)^{-y} D_0(r,x,y)\\
&=& \rho_t(m)^k \rho_t(\ell)^j D_t(r,x,y),
\end{eqnarray*}
so $D_t$ satisfies the right equivariance properties. Next, since $\rho_t(m)$ fixes $\Cln$ pointwise, we still have that
$\lim_{r \rightarrow 0} D_t(r,x,y)$ is independent of $x$. Further, for small $t$, $\lim_{r \rightarrow 0} D_t(r,x,y)$ will still be a local submersion to $\Cln$. So, $D_t$ is the developing map for a structure with cone-like singularities on a neighborhood of the singular locus $\Sigma$. Further, $D_t$ converges to $D_0$ in the $\mathscr C^1$ topology (in fact, in $\mathscr C^\infty$) on compacts of $\widetilde T$. Now, the definition of $D_t$ on $\widetilde T$ and the definition of $D_t$ on $\widetilde M_0$ may not agree on the overlap. So, we glue these two maps together using a bump function which is supported away from the singular locus. Finally, extend $D_t$ to the other lifts of $T$ in $\widetilde M$ by $\rho_t$ equivariance. This gives globally defined maps $D_t : \widetilde M \rightarrow \RP^3$ which converge in $\mathscr C^1$, on compacts, to $D_0$. Thus for sufficiently small $t$, the $D_t$ are local diffeomorphisms. 
\end{proof}

\subsection{Regeneration of $\mathbb H^3$ and $\AdS^3$ structures from $\HP^3$}
\label{HP-regen}

As the class of cone-like singularities specializes to cone singularities in the $\HH^3$ case, tachyons in the $\AdS^3$ case, and infinitesimal cone singularities in the $\HP^3$ case, we get the following regeneration statement immediately from Proposition~\ref{deform-cone-like}.

\begin{Proposition}[Regeneration with cone-like singularities]
Let $N$ be a closed three-manifold, with $\Sigma$ a knot, and let $M = N\setminus \Sigma$ with $m \in \pi_1 M$ the meridian around $\Sigma$.  Let $X$ be either $\XX_1 = \mathbb H^n$ or $\XX_{-1} = \AdS^n$. Let $\rho_t: \pi_1 M \rightarrow \Isom(X)$ be a family of representations defined for $t\geq 0$ such that 
\begin{itemize}
\item The path of conjugates $\resc_t \rho_t \resc_t^{-1}$ converges as $t  \rightarrow 0$ in $\mathscr C^1$ to a representation $\rho_{\HP}$,
\item $\rho_{\HP}$ is the holonomy of an $\HP$ structure on $N$ with infinitesimal cone singularity along $\Sigma$.
\item $\rho_t(m)$ is a rotation if $X = \HH^3$, or a boost if  $X = \AdS^3$.
\end{itemize}
Then, for sufficiently small $t > 0$, we can construct a family of $X$ structures on $N$ with singularities at $\Sigma$. For each $t$, the holonomy representation of the smooth part is $\rho_t$, and the structures have cone singularities if $X = \HH^3$ or tachyon singularities if $X = \AdS^3$.
\label{regen-sing}
\end{Proposition}
\begin{proof}
The proof is the same as the proof of Proposition~\ref{regen}. At time $t=0$ we have an $\HP$ structure with infinitesimal cone singularity. We regard this as a projective structure with cone-like singularities. If $X = \HH^3$, then $\sigma_t = \resc_t \rho_t \resc_t^{-1}$ is a representation landing in $G_t = \text{Isom}(\XX_t)$. The $\sigma_t$ limit to $\rho_{\HP}$. By Proposition~\ref{deform-cone-like} there is a family of cone-like projective structures very close to the $\HP$ structure that realize the $\sigma_t$ as holonomy (for short time). 
The developing maps of these structures map a compact fundamental domain $K$ (which includes the singularity) to a compact region inside of $\RP^3$ that for small $t$ is very close to the image of $K$ by the developing map of the $\HP$ structure. Thus, for all sufficiently small $t$, the image will lie inside of $\XX_t$ ensuring the developing maps define a family of $(\XX_t,G_t)$ structures. Applying the inverse of the rescaling map $\resc_t^{-1}$ gives a family of $\XX_1 = \HH^3$ structures with cone singularities. If $X = \AdS^3$ everything works the same, except that $\sigma_t$ lands in $G_{-t}$ and we get $(\XX_{-t}, G_{-t})$ structures that, by applying $\resc_t^{-1}$, are equivalent to $\AdS$ structures with tachyons.
\end{proof}

This proposition says that we can recover collapsing hyperbolic cone and $\AdS$ tachyon structures from an $\HP$ structure and a suitable path of representations. We use this proposition to prove our main regeneration theorem, Theorem~\ref{thm:main-regen} from the introduction, as follows:
We use the condition $H^1(\pi_1 M, \so(2,1)_{Ad \rho_0}) = \mathbb R$ to get representations into $\PSO(3,1)$ and $\PSO(2,2)$ satisfying the conditions of Proposition~\ref{regen-sing}. In the hyperbolic case, the proof of this makes use of the complex structure of the variety of $\PSO(3,1)$ representations coming from the isomorphism $\PSO(3,1) \cong \PSL(2,\CC)$. This isomorphism can be generalized to give $\PGL(2,\cdot)$ descriptions of all of the isometry groups $G_s$. Working with the $\PGL(2,\cdot)$ description of isometry groups allows for the most natural proof of Theorem~\ref{thm:main-regen}. So, we take a detour in the next section and give the proof of the Theorem in Section~\ref{proof-main-regen}.

\subsection{The $\PGL(2)$ description of isometry groups}
\label{PGL-description}

In dimension three, there is a useful alternative description of the isometry groups $G_t$ of our models $\XX_t$ (Section~\ref{projective-models}) which generalizes the isomorphism $\PSO(3,1) \cong \PSL(2,\CC)$.

Let $\BB_s = \mathbb R + \mathbb R \kappa_s$ be the real two-dimensional (commutative) algebra generated by a non-real element $\kappa_s$ with $\kappa_s^2 = -\text{sign}(s)s^2$. As a vector space $\BB_s$ is spanned by $1$ and $\kappa_s$. There is a conjugation action: 
$ \overline{(a + b \kappa_s)} := a - b \kappa_s  $
which defines a square-norm
$$ |a + b\kappa_s|^2 := (a + b\kappa_s)\overline{(a + b \kappa_s)} = a^2 -b^2 \kappa_s^2 \ \in \ \mathbb R.$$
Note that $|\cdot|^2$ may not be positive definite. We refer to $a$ as the \emph{real part}  and $b$ as the \emph{imaginary part} of $a + b\kappa_s$. It easy to check that $\BB_s$ is isomorphic to $\CC$ when $s > 0$.

In the case $s = -1$, we will denote $\kappa_s$ by the letter $\tau$.  It is also easy to check that when $s < 0$, $\BB_s$ is isomorphic to $\Rtau$.
In the case $s = 0$, we will denote $\kappa_s$ by the letter $\sigma$. The algebras $\Rtau$ and $\Rsigma$ play a central role in the study of ideal triangulations of $\AdS$ and $\HP$ manifolds, analogous to the role of the complex numbers in the study of ideal triangulations of hyperbolic manifolds. See \cite{Danciger-11}

Now consider the $2\times2$ matrices $M_2(\BB_s)$. 
Let $\Herm(2, \BB_s) = \{ A \in M_2(\BB_s) : A^* = A\}$ denote the $2\times 2$ Hermitian matrices, where $A^*$ is the conjugate transpose of $A$.  As a real vector space, $\Herm(2,\BB_s) \cong \mathbb R^4$. We define the following (real) inner product on $\Herm(2,\BB_s)$:
$$ \left\langle \minimatrix{a}{z}{\bar{z}}{d}, \minimatrix{e}{w}{\bar{w}}{h} \right\rangle = -\frac{1}{2}tr\left( \minimatrix{a}{z}{\bar{z}}{d} \minimatrix{h}{-w}{-\bar{w}}{e}\right).$$
The signature of this metric depends on $s$. 

\begin{Proposition}
The region $\XX_s$ in $\RP^3$ defined in Section~\ref{projective-models} can be alternately defined by
$$\XX_s = \left\{ X \in \Herm(2,\BB_s) : \langle X, X\rangle < 0 \right\} / X \sim \lambda X \text{ for } \lambda \in \mathbb R^*$$
where we use the coordinates $X = \minimatrix{x_1+x_2}{x_3 - x_4 \kappa_s}{x_3 + x_4\kappa_s}{x_1 - x_2}$ on $\Herm(2,\BB_s)$. Note that $\langle X, X\rangle = -\text{det}(X) = -x_1^2 + x_2^2 + x_3^2 -\kappa_s^2 x_4^2$.
\end{Proposition}

The ideal boundary $\partial^\infty \XX_s$, given by the projectivized light cone with respect to this metric, is exactly the projectivized rank one Hermitian matrices, where for a Hermitian matrix $X$, rank one means $\det(X) = 0$, $X \neq 0$. Any rank one Hermitian matrix $X$ can be decomposed (uniquely up to $\pm$) as $X = v v^*$ where $v \in \BB_s^2$ is a two-dimensional column vector with entries in $\BB_s$ (and $v^*$ denotes the transpose conjugate). Further $v$ must satisfy that $\lambda v = 0$ for $\lambda \in \BB_s$ if and only if $\lambda = 0$. This gives the identification $$\partial ^\infty \mathbb X = \mathbb P^1 \BB_s = \left\{ \btwovector{x}{y} : x\cdot \alpha  = 0 \text{ and }  y\cdot \alpha = 0 \text{ for } \alpha \in \BB_s \iff \alpha = 0 \right\} / \sim$$ where $\btwovector{x}{y} \sim \btwovector{x\lambda}{y \lambda} \text{ for } \lambda \in \BB_s^{\times}.$

\begin{Definition}
\label{defn-PGL+} We denote by $\PGL^+(2,\BB_s)$ the $2 \times 2$ matrices $A$ with entries in $\BB_s$ such that $|\det(A)|^2 > 0$, up to the equivalence $A \sim \lambda A$ for any $\lambda \in \BB_s^\times$.
\end{Definition}
Note that the condition $|\det(A)|^2 > 0$ is only needed in the case $s < 0$. For $s \geq 0$, $\PGL^+$ and $\PGL$ are the same. For $s > 0$, $\PGL^+$ is the same as $\PSL$.

We will think of $\PGL^+(2,\BB_s)$ as determinant $\pm 1$ matrices with entries in $\BB_s$ up to multiplication by a square root of $1$ (if $s < 0$, there will be four such square roots). 
We note that $\PGL^+(2, \BB_s)$ acts by Mobius transformations on $\partial^\infty \XX_s = \mathbb P^1 \BB_s$. This action extends to all of $\XX_s$, giving a map $\PGL^+(2,\BB_s) \rightarrow G_s = \text{Isom}(\XX_s)$, as follows:
$$ A \cdot X := A X A^* \text{ \ \ \ where $X \in \XX_s$ and $\det(A) =\pm 1$}.$$

\begin{Proposition}
For $s \neq 0$ The map $\PGL^+(2,\BB_s) \rightarrow G_s = \Isom^+(\XX_s)$ is an isomorphism. Note that in the case $s= 1$, this is the usual isomorphism $\PSL(2,\CC) \cong \PSO(3,1)$.
\end{Proposition}
\begin{proof}
The proof is an easy exercise. Use the coordinates $X = \minimatrix{x_1+x_2}{x_3 - x_4 \kappa_s}{x_3 + x_4\kappa_s}{x_1 - x_2}$ on $\Herm(2,\BB_s)$.
\end{proof}

\begin{Remark}
In fact, the orientation reversing isometries are also described by $\PGL^+(2, \BB)$ acting by $X \mapsto A \overline{X} A^*$.
\end{Remark}

Note that with the coordinates $X = \minimatrix{x_1+x_2}{x_3 - x_4 \kappa_s}{x_3 + x_4\kappa_s}{x_1 - x_2}$ on $\Herm(2,\BB_s)$, the rescaling map $\resc_s: \XX_1 \rightarrow \XX_s$ defined in Section~\ref{projective-models} corresponds to the algebraic rescaling map $\rescPGL_s: \CC = \BB_1 \rightarrow \BB_s$ defined by $i \mapsto \kappa_s/|s|$. This observation gives the following proposition:
\begin{Proposition}
\label{rescale-PSL}
For $s > 0$, $\rescPGL_s$ defines an isomorphism $\PSL(2,\BB_1) \rightarrow \PSL(2,\BB_s)$ which corresponds to the isomorphism $G_1 \rightarrow G_s$ given by conjugation by $\resc_s$. 
\begin{equation}
\label{PGL-O}
\xymatrix{\PSL(2,\CC) \ar[d]^{\cong} \ar[r]^{\rescPGL_s } & \PSL(2,\BB_s) \ar[d]^{\cong}\\ \PSO(3,1) \ar[r]^{\resc_s } & G_s}
\end{equation}
\end{Proposition}

Similarly, for $s < 0$, the rescaling map $\resc_s: \XX_{-1} \rightarrow \XX_s$ defined in Section~\ref{projective-models} corresponds to the algebraic rescaling map $\rescPGL_s:  \BB_{-1} \rightarrow \BB_s$ defined by $\tau  \mapsto \kappa_s/|s|$. Again, we get

\begin{Proposition}
\label{rescale-PSL}
For $s < 0$, $\rescPGL_s$ defines an isomorphism $\PGL^+(2,\BB_{-1}) \rightarrow \PGL^+(2,\BB_s)$, which corresponds to the isomorphism $G_1 \rightarrow G_s$ given by conjugation by $\resc_s$. 
\begin{equation}
\label{PGL-O}
\xymatrix{\PGL^+(2,\Rtau) \ar[d]^{\cong} \ar[r]^{\rescPGL_s } & \PGL^+(2,\BB_s) \ar[d]^{\cong}\\ \PSO(2,2) \ar[r]^{\resc_s } & G_s}
\end{equation}
\end{Proposition}

Recall that in the case $s=0$, the metric on $\XX_0$ is degenerate, so that the isometries of $\XX_0$ ended up being too large to be of use. The half-pipe group $G_{\HP}$ was defined to be a strict subgroup giving a useful structure for the purposes of geometric transitions.
\begin{Proposition}
The map $\PGL(2,\Rsigma) \rightarrow G_0 = \text{Isom}(\XX_0)$ maps $\PGL(2,\Rsigma)$ isomorphically onto $G_{\HP}^+$.
\end{Proposition}
\begin{proof}
To begin, we think of $\Rsigma$ as the cotangent bundle of $\RR$. The element $\sigma$ should be thought of as an infinitesimal quantity, whose square is zero. Similarly, $\PGL(2,\Rsigma)$ is the cotangent bundle of $\PGL(2,\RR)$:
\begin{Lemma}
Let $A + B\sigma$ have determinant $\pm 1$. Then $\det A = \det(A + B\sigma) = \pm 1$ and $\operatorname{tr} A^{-1} B = 0$. In other words $B$ is in the tangent space at $A$ of the matrices of constant determinant $\pm 1$.
\end{Lemma}
Any element of $\Herm(2,\Rsigma)$ can be expressed uniquely in the form $X + Y\sigma$ where $X = \minimatrix{x_1+x_2}{x_3}{x_3}{x_1-x_2} = X^T$ is symmetric and $Y = \minimatrix{0}{-x_4}{x_4}{0} = -Y^T$ is skew-symmetric. Then
\begin{eqnarray*}
(A+B\sigma)(X+Y\sigma)(A + B \sigma)^* &=& (A+B\sigma)(X+Y\sigma)(A^T - B^T \sigma) \\
&=& AXA^T + (BXA^T - AXB^T +AYA^T)\sigma 
\end{eqnarray*}
where we note that $AXA^T$ is symmetric and $(BXA^T - AXB^T +AYA^T)$ is skew-symmetric.
The symmetric part $X \mapsto AXA^{T}$, written in coordinates gives the familiar isomorphism $\Phi : \PGL(2, \RR) \rightarrow \OO(2,1)$. In $(x_1,x_2,x_3,x_4)$ coordinates the transformation defined by $A + B\sigma$ has matrix
\begin{equation*}
\begin{pmatrix}
\Phi(A) & 0\\ v(A,B) & c(A,B)
\end{pmatrix}
\end{equation*}
The skew-symmetric part $X+Y\sigma \mapsto (BXA^T - AXB^T +AYA^T)\sigma$, written in $(x_1,x_2,x_3,x_4)$ coordinates gives the bottom row of this matrix:
\begin{eqnarray*} v(A,B) &=& \begin{pmatrix} (-ce - df + ag +bh) & (-ce + df + ag -bh) & (-cf - de + ah +bg)  \end{pmatrix} \\ c(A,B) &=& \det(A) \ =\  \det(\Phi(A)) \ =\  \pm 1
\end{eqnarray*}
where $A = \minimatrix{a}{b}{c}{d}$, $B = \minimatrix{e}{f}{g}{h}$. To show that $\PGL(2,\Rsigma) \rightarrow G_{\HP}^+$ is an isomorphism, one easily checks that given $A$, the map $B \rightarrow v(A,B)$ is a linear isomorphism to $\RR^3$.
\end{proof}

Finally, we restate the condition of \emph{compatibility to first order} (Equation~\ref{limit-hol} of Section~\ref{rescaling-defnHP}) in these terms. 
In order to make sense of continuity and limits for paths of representations over the varying algebras $\BB_s$, we can embed all of the $\BB_s$ in a larger Clifford algebra $\mathcal C$ (see \cite{Danciger-11}, Section~5.4). For our purposes here, we assume that $\kappa_s \rightarrow \kappa_0$ as $s \rightarrow 0$.

In the hyperbolic case:
\begin{Proposition}
\label{compatibility-PGLC}
Let $\rho_t: \pi_1 M \rightarrow G_{+1}$ be a path of representations, defined for $t \geq 0$, converging to a representation $\rho_0$ with image in the subgroup $H_0 = \begin{pmatrix} \OO(2,1) & 0 \\ 0 & \pm 1 \end{pmatrix}$. Then the corresponding representations $\tilde \rho_t : \pi_1 M \rightarrow \PGL(2,\CC)$ limit to a representation $\tilde \rho_0$ into $\PGL(2,\RR)$. Suppose further that $\resc_t \rho_t \resc_t^{-1}$ limit to a representation $\rho_{\HP}.$ Then $\rescPGL_t \tilde \rho_t \rescPGL_t^{-1} \xrightarrow[t\rightarrow 0]{} \tilde \rho_{\HP}$ where $\tilde \rho_{\HP}$ is the representation into $\PGL(2,\Rsigma)$ corresponding to $\rho_{\HP}$. Further $\tilde \rho_{\HP}$ is defined by
\begin{align} \label{eqn:rhoHPfromC}
\tilde \rho_{\HP} &= \tilde \rho_0 + \sigma \frac{d}{dt} \operatorname{Im} \tilde \rho_t \Big |_{t=0}.
\end{align}
\end{Proposition}

Similarly, in the $\AdS$ case:
\begin{Proposition}
\label{compatibility-PGLRtau}
Let $\rho_t: \pi_1 M \rightarrow G_{-1}$ be a path of representations, defined for $ t\leq0$, converging to a representation $\rho_0$ with image in the subgroup $H_0 = \begin{pmatrix} \OO(2,1) & 0 \\ 0 & \pm 1 \end{pmatrix}$. Then the corresponding representations $\tilde \rho_t : \pi_1 M \rightarrow \PGL(2,\Rtau)$ limit to a representation $\tilde \rho_0$ into $\PGL(2,\RR)$. Suppose further that $\resc_t \rho_t \resc_t^{-1}$ limit to a representation $\rho_{\HP}.$ Then $\rescPGL_t \tilde \rho_t \rescPGL_t^{-1} \xrightarrow[t\rightarrow 0]{} \tilde \rho_{\HP}$ where $\tilde \rho_{\HP}$ is the representation into $\PGL(2,\Rsigma)$ corresponding to $\rho_{\HP}$. Further $\tilde \rho_{\HP}$ is defined by
\begin{align}\label{eqn:rhoHPfromRtau}
\tilde \rho_{\HP} &= \tilde \rho_0 + \sigma \frac{d}{dt} \operatorname{Im} \tilde \rho_t \Big |_{t=0}.
\end{align}
\end{Proposition}

In both propositions, Equations~\ref{eqn:rhoHPfromC} and \ref{eqn:rhoHPfromRtau} are made sense of by choosing (for each $\gamma \in \pi_1M$) a lift of $\tilde \rho_t$ to $\GL(2, \cdot)$ with constant determinant $\pm 1$ such that the limit $\tilde \rho_0 \in \GL(2,\RR)$ is real.


\subsection{Proof of regeneration theorem}
\label{proof-main-regen}

We restate Theorem~\ref{thm:main-regen} in terms of $\PGL(2,\BB_s)$ isometry groups.

\begin{Theorem}
Let $N$ be a closed orientable $\HP^3$ manifold with infinitesimal cone singularity of infinitesimal angle $-\omega$ along the knot $\Sigma$. Let $M = N\setminus \Sigma$ be the smooth part and let $\rho_{\HP}: \pi_1 M \rightarrow \PGL(2,\Rsigma)$ be the holonomy representation. Suppose that the real part $\rho_0$ of $\rho_{\HP}$ satisfies the condition $H^1(\pi_1 M, \ssl(2,\RR)_{Ad \rho_0}) = \mathbb R$. Then there exists singular geometric structures on $(N,\Sigma)$ parametrized by $t \in (-\delta, \delta)$ which are
\begin{itemize}
\item hyperbolic cone structures with cone angle $2\pi - \omega t$ for $t > 0$
\item $\AdS$ structures with a tachyon of mass $-\omega t$ for $t < 0$.
\end{itemize}
\end{Theorem}

\begin{proof}

We begin with a lemma about the representation variety $\mathscr R(\pi_1 M, \PGL(2,\RR))$ of representations modulo conjugation.
\begin{Lemma}
The condition $H^1(\pi_1 M, \ssl(2,\RR)_{Ad \rho_0}) = \mathbb R$ guarantees that the representation variety $\mathscr R(\pi_1 M, \PGL(2,\RR))$ is smooth at $\rho_0$.
\end{Lemma}
\begin{proof}
It is a standard fact that $H^1(\pi_1 M, \ssl(2,\RR)_{Ad \rho}) \rightarrow H^1(\pi_1 \partial M,  \ssl(2,\RR)_{Ad \rho})$ has half-dimensional image (see for example \cite{Hodgson-05}). In this case, $\rho_0(m) = 1$ and $\rho_0(\ell)$ is a non-trivial translation (possibly plus a rotation by $\pi$), so any nearby representation $\phi$ of $\pi_1(\partial M)$ preserves an axis so that $H^1(\pi_1 \partial M,  \ssl(2,\RR)_{Ad \phi})$ has dimension equal to two. So, $H^1(\pi_1 M, \ssl(2,\RR)_{Ad \rho})$ has dimension at least one for all $\rho$ nearby $\rho_0$. It follows that $\rho_0$ is a smooth point of $\mathscr R(\pi_1 M, \PGL(2,\RR))$, and that the tangent space at $\rho_0$ is one dimensional. 
\end{proof} 

Let $m$ be a meridian around $\Sigma$ in the direction consistent with the orientation of $\Sigma$ (so that in particular, the discrete rotational part of the holonomy of $m$ is $+2\pi$).

\vspace{0.1in}
{\noindent \bf Hyperbolic case ($t > 0$):}
 In order to use Proposition~\ref{regen-sing}, we must produce for $t > 0$ a path of representations $\rho_t$ into $\PSL(2,\CC)$ with the following properties:
\begin{enumerate}
\item $\rho_t \rightarrow \rho_0$
\item $\rho_t(m)$ is a rotation by $2\pi-\omega t$ \label{prop1}
\item $\resc_t \rho_t \resc_t^{-1}$ converges to $\rho_{\HP}$ as $t \rightarrow 0$. By Proposition~\ref{compatibility-PGLC} this is equivalent to $$\frac{d}{dt} \operatorname{Im} \rho_t \Big |_{t=0} = \operatorname{Im} \rho_{\HP}.$$
\end{enumerate}

Now, our $\HP$ representation gives a $\PGL(2,\RR)$ tangent vector at $\rho_0$ as follows:
$\rho_{\HP}(\gamma) = \rho_0(\gamma) + Y(\gamma)\ \sigma.$ Define $z(\gamma) = Y(\gamma)\rho_0(\gamma)^{-1}$. Then $z$ is an $\ssl(2,\RR)_{Ad \rho_0}$ co-cycle. It spans the tangent space of $\mathscr R(\pi_1 M, \PGL(2,\RR))$.  As the structure is singular, we must have $z(m) \neq 0$. Thus the translation length of $m$ increases (or decreases) away from zero. The complexified variety $\mathscr R(\pi_1 M, \PSL(2,\CC))$ is also smooth at $\rho_0$ and $\mathscr R(\pi_1 M, \PSL(2,\CC)) \rightarrow \mathscr R(\pi_1 \partial M, \PSL(2,\CC))$ is a local immersion at $\rho_0$. The variety $\mathscr R(\pi_1 \partial M, \PSL(2,\CC))$ has complex dimension $2$.
\begin{Lemma}
The subset $$S = \{ \rho \in \mathscr R(\pi_1 \partial M, \PSL(2,\CC)) : \rho(m) \text{ is elliptic}\}$$ is locally a smooth real sub-manifold of dimension three. The image of $\mathscr R(\pi_1 M, \PSL(2,\CC))$ in $\mathscr R(\pi_1 \partial M, \PSL(2,\CC))$ is transverse to $S$.
\end{Lemma}
\begin{proof}[Proof of lemma]
That $S$ is smooth of dimension three follows immediately from the fact that $\mathscr R(\pi_1 \partial M, \PSL(2,\CC))$ is parameterized (near $\rho_0$) by the complex lengths of $m,\ell$. The image of $\mathscr R(\pi_1 M, \PSL(2,\CC))$ in $\mathscr R(\pi_1 \partial M, \PSL(2,\CC))$ is transverse to $S$ because $z$ increases translation length of $m$ away from zero.
\end{proof}
Now, from the lemma, we have that the $\PSL(2,\CC)$ representations of $\pi_1 M$ with $m$ elliptic near $\rho_0$ form a smooth real one-dimensional manifold. The tangent space at $\rho_0$ is spanned by $iz(\cdot)$. Thus the rotation angle of $m$ is changing along this manifold and we can choose $\rho_t$ as desired. 

\vspace{0.1in}

{\noindent \bf $\AdS$ case ($t < 0$):}
We obtain, from the argument above, a path $\phi_t: \pi_1 M \rightarrow \PGL(2,\RR)$ defined in a neighborhood of $t=0$ with $\phi_0 = \rho_0$, $\frac{d}{dt}\phi_t \big |_{t=0} = \text{Im} \rho_{\HP}$ and $z(m) =\frac{d}{dt}\phi_t(m)$ is an infinitesimal translation by $-\omega$ along the axis $\Cln$ of $\rho_0(\ell)$. 
We may assume that the axis in $\HH^2$ preserved by $\phi_t(\partial M)$ is also $\Cln$ (constant). Now, define $\rho_t : \pi_1 M \rightarrow \PGL^+(2,\Rtau)$ by 
\begin{equation*}
\rho_t(\cdot) = \frac{1+\tau}{2}\phi_t(\cdot) + \frac{1-\tau}{2} \phi_{-t}(\cdot).
\end{equation*}
A quick computation shows that $\frac{d}{dt} \rho_t \big |_{t=0} =  \tau \frac{d}{dt}\phi_t \big |_{t=0}$. Further, $\rho_t(m)$ is a boost around the axis $\Cln$ by hyperbolic angle $\omega t$. So Proposition~\ref{regen-sing} implies the result for $t< 0$.

\end{proof}

%
%

\section{Hyperbolic and AdS manifolds that collapse onto the $(2,m,m)$ triangle orbifold}
\label{example1}

In this section we prove Theorem~\ref{thm:regen-2mm}, restated here for convenience:

\begin{reptheorem}{thm:regen-2mm}
Let $m \geq 5$, and let $S$ be the hyperbolic structure on the two-sphere with three cone points of order $2,m,m$.  Let $N$ be the unit tangent bundle of $S$. Then, there exists a knot $\Sigma$ in $N$ and a continuous path of real projective structures $\mathscr P_t$ on $N$, singular along $\Sigma$, such that $\mathscr P_t$ is conjugate to:
\begin{itemize}
\item
a hyperbolic cone structure of cone angle $\alpha < 2\pi$, when $t > 0$, or
\item
an $\AdS$ structure with tachyon singularity of mass $\phi < 0$, when $t < 0$.
\end{itemize}
As $t \to 0$, the cone angle $\alpha \to 2 \pi^-$ (resp. $\phi \to 0^-$) and the hyperbolic geometry (resp. $\AdS$ geometry) of $\mathscr P_t$ collapses to the surface $S$. 
\end{reptheorem}

There are two steps to prove the theorem. First, we construct a half-pipe structure on $N$ with infinitesimal cone singularity along a circle $\Sigma$. This construction is explicit and gives a first order approximation of the geometry of the nearly collapsed manifolds in the Theorem. Second, we give a simple analysis of the representation variety and apply Theorem~\ref{thm:main-regen}.

The methods of this section should be compared/contrasted with \cite{Porti-10}, in which hyperbolic cone orbifolds that collapse to a hyperbolic polygon are shown to exist.

\subsection{Half-pipe geometry in dimension three}
\label{s:HP3-geometry}
Before getting to the main construction of this section, we give some lemmas useful for working with $\HP$ geometry in dimension three.
It will be most convenient to work with the model given in Section~\ref{PGL-description} (with $s=0$). Recall the algebra $\RR + \RR\sigma$, with $\sigma^2 = 0$. The half-pipe space is given by 
$$\HP^3:= \mathbb X_0 = \left\{ X + Y \sigma : X,Y \in M_2(\RR), X^T = X, \det(X) > 0, Y^T = -Y\right\}/\sim $$
where $(X + Y \sigma) \sim \lambda (X + Y \sigma)$ for $\lambda \in \RR^\times$. The diffeomorphism $\mathbb X_0 \to \HH^2 \times \RR$ of Section~\ref{HP-geometry} is given in these coordinates by 
\begin{align}\label{product-coords-dim3}
X + Y \sigma \mapsto (X, L)
\end{align} where we interpret the symmetric matrices $X$ of positive determinant, up to scale, as a copy of $\HH^2$ in the usual way and the length $L$ along the fiber is defined by the equation 
\begin{align} \label{fiber-length-dim3}
Y &= L \sqrt{\det X}\begin{pmatrix} 0 & -1\\ 1 & 0 \end{pmatrix}. 
\end{align}
We restrict to the identity component of the structure group, which is given by
\begin{align*}
G_0 &= \PSL(2, \RR + \RR\sigma)\\
&= \{ A + B \sigma : A \in \SL(2,\RR), \text{ and } B \in T_A \SL(2,\RR) \} / \pm.
\end{align*}
The structure group identifies with the tangent bundle $T \PSL(2,\RR)$, and it will be convenient to think of elements as having a finite component $A\in \PSL(2,\RR)$ and an infinitesimal component $a \in \mathfrak{sl}(2,\RR)$, via the isomorphism
\begin{align*}
\PSL(2,\RR) \ltimes \mathfrak{sl}(2,\RR) &\to G_0\\
(A,a) &\mapsto A + Aa \sigma
\end{align*}
where $Aa \in T_A \PSL(2,\RR)$.
(This is the usual isomorphism $G \ltimes \mathfrak g \to T G$ for a Lie group $G$ with Lie algebra $\mathfrak g = T_1 G$, gotten by left translating vectors from the identity.)
Thinking of $a \in \mathfrak{sl}(2,\RR)$ as an infinitesimal isometry of $\HH^2$, recall that at each point $X \in \HH^2$ we may decompose $a$ into its translational ($X$-symmetric) and rotational ($X$-skew) parts: 
\begin{align*}
a &= a_{X\text{-sym}} + a_{X\text{-skew}} \\ 
&:= \frac{1}{2}\left( a + Xa^T X^{-1} \right) + \frac{1}{2} \left( a - Xa^T X^{-1} \right).
\end{align*}
where the rotational part $a_{X\text{-skew}}$ is a rotation centered at $X$ of infinitesimal angle $\rot(a,X)$ defined by
$$\sqrt{X}^{-1} a_{X\text{-skew}} \sqrt{X} = \rot(a,X) \begin{pmatrix} 0 & -1/2 \\ 1/2 & 0\end{pmatrix},$$
 The action of an element of $G_0$ in the fiber direction depends on the rotational part of the infinitesimal part of that element.
\begin{figure}[h]
{\centering

\def\svgwidth{3.7in}
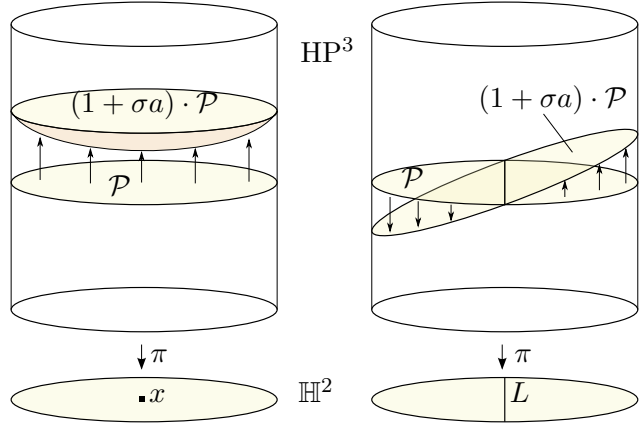

\caption{The action of $1+a\sigma$ on $\HP^3$ when $a$ is an infinitesimal rotation centered at $x$ (left) or $a$ is an infinitesimal translation along $L$ (right).}
}
\end{figure}

\begin{Lemma}
\label{lem:HP3-action}
The action of a pure infinitesimal $1 + a\sigma$ on the point $X + Y\sigma \in \mathbb X_0$ is by translation in the fiber direction by amount equal to the rotational part $\rot(a, X)$ of the infinitesimal isometry $a$ at the point $X \in \HH^2$.  In the product coordinates~(\ref{product-coords-dim3}):
$$1+ a \sigma : (X,L) \mapsto (X, L + \rot(a,X)).$$
More generally, the action of $A + Aa\sigma$ is given by 
$$A + Aa\sigma : (X,L)\mapsto (A\cdot X, L + \rot(a,X)).$$
\end{Lemma}

\begin{proof}
\begin{align*}
(1+a\sigma)\cdot (X + Y\sigma) &= (1+a\sigma) (X+  \sigma Y) (1 - a^T\sigma)\\
&= X + \sigma Y + \sigma(a X - X a^T)\\
&= X + \sigma Y + \sigma \ 2 a_{X\text{-skew}} X\\
&= X + \sigma Y + \sigma \ 2 \rot(a,X) \sqrt{X} \begin{pmatrix} 0 & -1/2\\ 1/2 & 0\end{pmatrix} \sqrt{X}\\
&= X + \sigma Y + \sigma \  \rot(a,X) \det(\sqrt{X}) \begin{pmatrix} 0 & -1\\ 1 & 0\end{pmatrix},
\end{align*}
and the first statement now follows from Equation~(\ref{fiber-length-dim3}). The second more general formula follows easily after left multiplication by $A$.
\end{proof}

Finally, we mention two very easy, but useful facts:
\begin{Lemma}
\label{lem:inf-commute}
Infinitesimals in $G_0$ commute: $$(1+a\sigma)(1+b\sigma)  = (1+b\sigma)(1+a\sigma).$$
\end{Lemma}
\begin{Lemma}
\label{lem:kHkinv}
Conjugating an infinitesimal $1+h \sigma$ by $K + Kk\sigma$ only depends on the finite part $K$ of $K+Kk\sigma$:
$$(K+Kk \sigma) (1+h\sigma) (K+Kk\sigma)^{-1} = K(1+h\sigma)K^{-1} = 1 + Ad_K h \ \sigma.$$
\end{Lemma}

\subsection{Building blocks for the unit tangent bundle $N$}
We now begin our half-pipe geometry construction of the unit tangent bundle $N$ of $S$.
\begin{figure}[h]
{\centering

\def\svgwidth{3.2in}
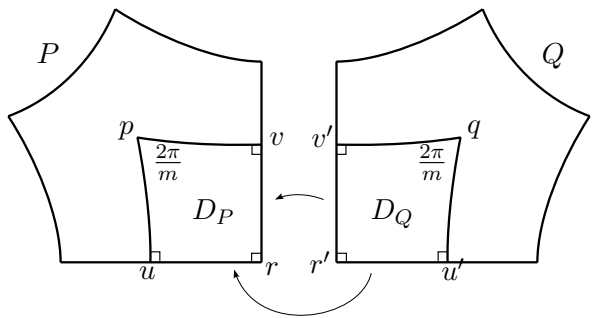

\caption{$S$ is constructed by glueing together two pieces $S_P$ and $S_Q$, which are quotients of the regular right-angled $m$-gons $P$ and $Q$ respectively. $D_P$ and $D_Q$ are fundamental domains for the rotational symmetry.\label{fig:DPDQ}}
}
\end{figure}
The hyperbolic orbifold surface $S$ can be constructed as follows. Take two isometric regular, right-angled $m$-gons, $P$ and $Q$. Let $A, B \in \PSL(2,\RR)$ be the counter-clockwise rotations of order $m$ about the center points $p$ and $q$ of $P$ and $Q$ respectively. Then $S_P = P/\langle A\rangle$ and $S_Q = Q/ \langle B \rangle$ are isometric monogons, each with an order $m$ interior cone point, and each having a right angle at its exterior vertex. Glueing $S_P$ and $S_Q$ together gives the $(2,m,m)$ triangle orbifold $S$. We choose a fundamental domain $D_P$ for $S_P$ as follows. Let $r$ be a corner point of $P$ and let $u$ and $v$ be the midpoints of the two edges adjacent to $r$. Then we define $D_P$ to be the quadrilateral spanned by $p,u,r,v$ as in Figure~\ref{fig:DPDQ}. Similarly, define the fundamental domain $D_Q$ to be the quadrilateral with vertices $q,u',r',v'$ as in the figure.  Now, $S$ is constructed by glueing $D_P$ to itself via $A$, glueing $D_Q$ to itself via $B$ and glueing $D_P$ to $D_Q$, identifying the edge $r'v'$ to $rv$ and the edge $r'u'$ to $ru$ as in the figure (using, in each case, the unique orientation preserving isometry).

\begin{figure}[h]
{\centering

\def\svgwidth{4.0in}
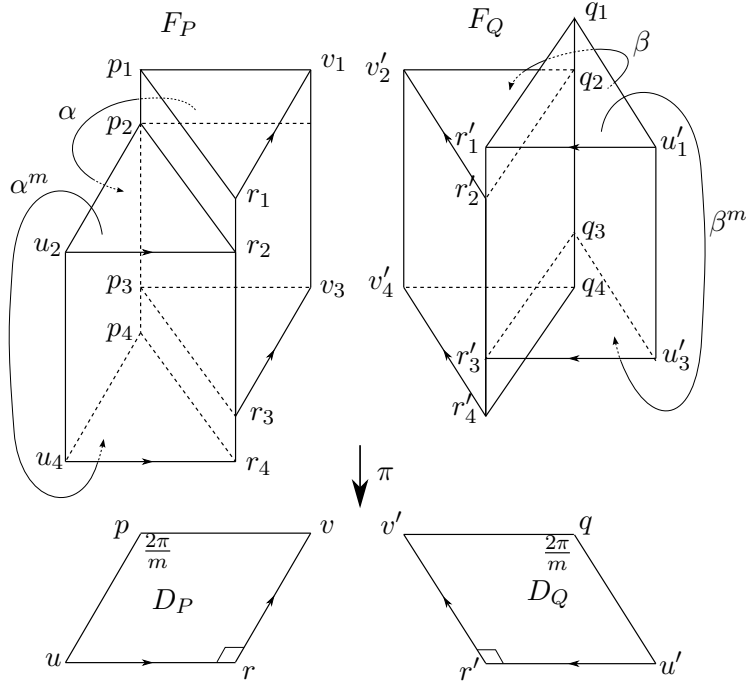

\caption{Build $N$ out of two pieces $F_P$ and $F_Q$. \label{2mm-FPFQ}}
}
\end{figure}

Next, we construct the unit tangent bundle $N$ of $S$ out of two pieces in a similar way. These pieces are suitable polyhedral lifts $F_P$ and $F_Q$ of $D_P$ and $D_Q$ to $\HP^3$. First we define $F_P$. Let $\alpha = A + A a \sigma$ be a lift of $A$ to $\PSL(2,\RR + \RR\sigma)$. Thinking of $\alpha$ as an element of $T \PSL(2,\RR)$ (via the discussion in Section~\ref{s:HP3-geometry}), the infinitesimal part $Aa \sigma$ of $\alpha$ describes the change $\dot \theta_\alpha$ in the rotation angle of $A$ (via the linearization of the usual trace relationship). The element $\alpha^m$ is an infinitesimal rotation of infinitesimal angle $m \dot \theta_\alpha$ centered at $p$ (the fixed point of $A$). Then, by Lemma~\ref{lem:HP3-action}, the action of $\alpha^m$ on $\HP^3$ is a translation in the fiber above each point $x\in \HH^2$ by amount
$$m \dot \theta_\alpha \cosh d(p,x)$$
where $d(\cdot , \cdot )$ is the distance function. We assume henceforth that $\dot \theta_\alpha > 0$ so that $\alpha^m$ translates upward in all fibers. 
Now refer to Figure~\ref{2mm-FPFQ} for the construction of $F_P$.
Let the lifts $ p_1,  p_2,  p_3,  p_4$ of $p$, the lifts $ r_1,  r_2,  r_3,  r_4$ of $r$, the lifts $ v_1,  v_3$ of $v$, and the lifts $ u_2,  u_4$ of $u$ satisfy the following:
\begin{align*}
 p_3 &= \alpha p_4 	&  p_1 = \alpha  p_2 \\
 v_3 &= \alpha u_4 	&  v_1 = \alpha  u_2 \\
 p_2 &= \alpha^m  p_4	&  p_1 = \alpha^m  p_3\\
 r_2 &= \alpha^m  r_4	&  r_1 = \alpha^m  r_3\\
 u_2 &= \alpha^m  u_4	&  v_1 = \alpha^m  v_3\\ 
\alpha r_4, r_3, v_3 & \text{ are colinear.} & \alpha r_2, r_1, v_1 \text{ are colinear.}
\end{align*}
All vertices of $F_P$ are determined by, say $ u_4,  r_4,  p_4$. The bottom of $F_P$, consisting of the hyperbolic triangles $\Delta  p_4  u_4  r_4$, $\Delta  p_3  r_3  v_3$, and the vertical ($\HP^2$) quadrilateral $\square p_4 r_4  r_3  p_3$, is glued by $\alpha^m$ to the top of $F_P$, consisting of the hyperbolic triangles $\Delta  p_2  u_2  r_2$, $\Delta  p_1  r_1  v_1$, and the vertical ($\HP^2$) quadrilateral $\square p_2 r_2 r_1  p_1$. The vertical ($\HP^2$) side face $\square p_4 u_4 u_2  p_2$ is glued by $\alpha$ to the adjacent vertical face $\square p_3  v_3  v_1  p_1$. The result is an $\HP$ structure on the unit tangent bundle $N_P$ of $S_P$. Note that the specific geometry of the bottom and top faces is not important (provided that the top is related to the bottom by $\alpha^m$).

Similarly, to construct $F_Q$, choose a lift $\beta = B + Bb\sigma$ of $B$, such that the infinitesimal change in rotation angle $\dot \theta_\beta > 0$. Let the lifts $ q_1,  q_2,  q_3,  q_4$ of $q$, $ r_1',  r_2',  r_3',  r_4'$ of $r'$, $ v_2',  v_4'$ of $v'$, and $ u_1',  u_3'$ of $u'$ satisfy the following:

\begin{align*}
 q_3 &= \beta q_4 	&  q_1 = \beta  q_2 \\
 u_3' &= \beta v_4' 	&  u_1' = \beta  v_2' \\
 q_2 &= \beta^m  q_4	&  q_1 = \beta^m  q_3\\
 r_2' &= \beta^m  r_4'	&  r_1' = \beta^m  r_3'\\
 u_1' &= \beta^m  u_3'	&  v_2' = \beta^m  v_4'\\ 
\beta r_4', r_3', u_3' & \text{ are colinear.} & \beta r_2', r_1', u_1' \text{ are colinear.}
\end{align*}
All vertices of $F_Q$ are determined by, say, $ v_4',  r_4',  q_4$. Now, the bottom of $F_Q$, consisting of the hyperbolic triangles $\Delta  q_4  v_4'  r_4'$, $\Delta  q_3  r_3'  u_3'$, and the vertical ($\HP^2$) quadrilateral $\square q_4 r_4'  r_3'  q_3$, is glued by $\beta^m$ to the top of $F_P$, consisting of the hyperbolic triangles $\Delta  q_2  v_2'  r_2'$, $\Delta  q_1  r_1'  u_1'$, and the vertical ($\HP^2$) quadrilateral $\square q_2 r_2' r_1'  q_1$. The vertical ($\HP^2$) side face $\square q_4 v_4' v_2'  q_2$ is glued by $\beta$ to the adjacent vertical face $\square q_3  u_3'  u_1'  q_1$. The result is an $\HP$ structure on the unit tangent bundle $N_Q$ of $S_Q$.

Next, we attempt to glue together the $\HP$ structures $N_P$ and $N_Q$ to produce an $\HP$ structure on $N$. First of all, we need to glue $\square r_2' r_4' v_4' v_2'$ to $\square r_1 r_3 v_3 v_1$. So in particular, we need the fiber of $F_P$ above a point $x$ of edge $r v$ to have the same length as the fiber of $F_Q$ above the corresponding point $y$ of edge $r' v'$ under the glueing
In other words, the action of $\alpha^m$ must translate the fiber $\pi^{-1}x$ the same amount that the action of $\beta^m$ translates the fiber $\pi^{-1} y$. This leads to the equation

\begin{align*}
m \dot \theta_\alpha \cosh d(x,p) &= m \dot \theta_\beta \cosh d(y,q)\\
\implies \dot \theta_\alpha &= \dot \theta_\beta
\end{align*}
where the second equality follows because $d(p,x) = d(q,y)$ by symmetry. In fact, this condition also ensures that $\square r_1' r_3' v_3' v_1'$ can be glued to $\square r_2 r_4 v_4 v_2$.
However, even if this is the case, glueing $N_P$ to $N_Q$ can not give a smooth structure. For let $g \in\PSL(\RR + \RR\sigma)$ map $\square r_2' r_4' v_4' v_2'$ to $\square r_1 r_3 v_3 v_1$. Note that $g$ is a lift of the map glueing edge $r'v'$ of $D_Q$ to edge $rv$ of $D_P$. Then $g \beta^m g^{-1}\alpha^{-m}$ is non-trivial as $g \beta^m g^{-1}$ is an infinitesimal rotation around $\pi_*(g) q \neq p$ (where $\pi_*(g)$ is the finite part of $g$). Hence, the corresponding closed curve $\gamma$ in $N$ which runs along a generic fiber of $N_P$, then runs once backwards around a generic fiber of $N_Q$ (passing between $N_P$ and $N_Q$ via $g$) has non-trivial holonomy. Since $\gamma$ is trivial in $\pi_1 N$, this glueing can not define a smooth $\HP$ structure on $N$; the curve $\gamma$ encircles an infinitesimal cone singularity in the face $\square r_1 r_3 v_3 v_1$.

\subsubsection{The singular $\HP$ structure on $N$: defining equations}

\begin{figure}[h]
{\centering

\def\svgwidth{4.5in}
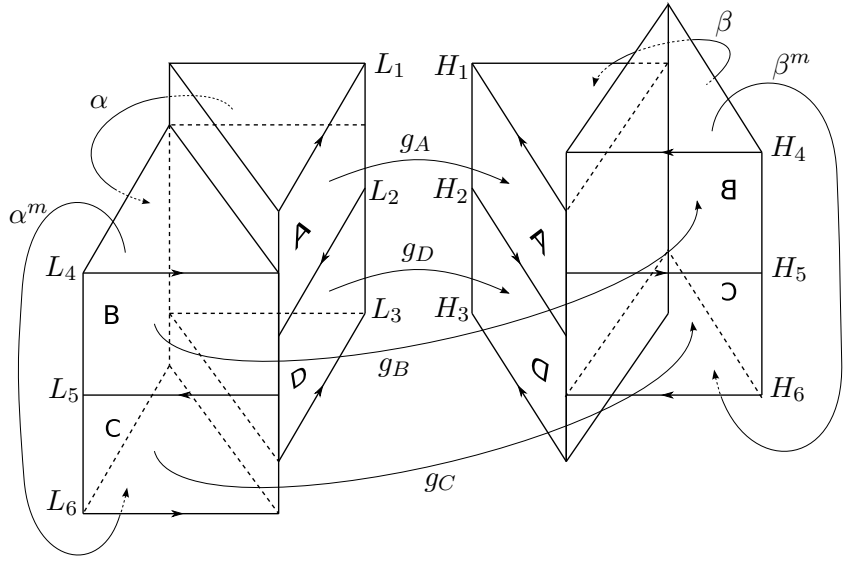

\caption{Glueing diagram for  the singular structure on $N$. The arrows represent paths crossing through a given face in the manifold. Each arrow is labeled with the corresponding glueing map. Note that glueing maps are dual to paths, e.g. the path from $F_P$ to $F_Q$ crossing though face $A$ is dual to the glueing map $g_A$ which glues $F_Q$ to $F_P$ along $A$. \label{2mmglueing1}}
}
\end{figure}
\begin{figure}[h]
{\centering

\def\svgwidth{2.3in}
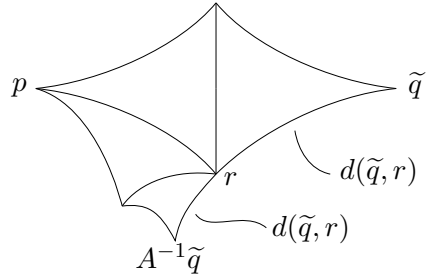

\caption{$r$ is the midpoint of the segment connecting $\widetilde q$ to $A^{-1} \widetilde q$.\label{2mmglueing2}}
}
\end{figure}

We add an edge to each of the faces $\square r_2' r_4' v_4' v_2'$, $\square r_1 r_3 v_3 v_1$, $\square r_2' r_4' v_4' v_2'$ and $\square r_1 r_3 v_3 v_1$ so that each quadrilateral is divided into two quadrilaterals stacked vertically. 
There are four glueing maps to be determined, $g_A, g_B, g_C, g_D$ which glue $F_Q$ to $F_P$ along $A, B, C, D$ respectively as in Figure~\ref{2mmglueing2}. The directed edges in the figure unite after glueing to form a simple closed (and oriented) curve $\Sigma$ in $N$. The goal here is to construct an $\HP$ structure on $N$ with infinitesimal cone singularity along $\Sigma$. The first condition we impose is that $\Sigma$ be totally geodesic (and smooth) in the glued up manifold. In the following equations, which we use in the next subsection, $L_i$ refers to the line containing a given line segment as labeled in Figure~\ref{2mmglueing1}.
%
\begin{align*}
 L_1 &= g_A H_1 &L_2 &= g_A H_2 & L_4 &= g_A H_5 \\
 L_2 &= g_D H_2 & L_3 &= g_D H_3 & L_5 &=g_D H_6\\
 L_4 &= g_B H_4 & L_5 &= g_B H_5 & L_2 &= g_B H_1\\
 L_5 &= g_C H_5 & L_6 &= g_C H_6 & L_3 &= g_C H_2\\
\alpha L_4 &= L_1 & \alpha L_5 &= L_2 & \alpha L_6 &= L_3\\
\beta H_1 &= H_4 & \beta H_2 &= H_5 & \beta H_3 &= H_6\\
\alpha^m L_3 &= L_1 & \alpha^m L_6 &= L_4\\
\beta^m H_3 &= H_1 & \beta^m H_6 &= H_4
\end{align*}

Next, we must pay special attention to the fiber $\mathcal F$ above $r$ in $F_P$ (the vertical edge at the intersection of face $A$ and $B$ or $D$ and $C$ in $F_P$). This is the exceptional fiber of order $2$. We must check that the $\HP$ structure is smooth at points of this fiber not belonging to one of the $L_i$. The following equations express that the holonomy of any curve encircling such a piece of $\mathcal F$ is trivial. In other words, the order $2$ exceptional fiber is non-singular:

\begin{equation}
\label{eqn:order2-fiber}
\begin{aligned}
g_A g_C^{-1} g_D \beta^{-m} g_B^{-1} &= 1	\hspace{0.5in}	& g_D g_C^{-1} \alpha^{-m} g_A g_B^{-1} &= 1\\
g_D \beta^{-m} g_B^{-1} g_A g_C^{-1} &=1		& g_A g_B^{-1} g_D g_C^{-1} \alpha^{-m} &= 1
\end{aligned}
\end{equation}

Next, (segments of) the $L_i$ become the singular locus $\Sigma$ after glueing, which has one component. Therefore all curves encircling pieces of the $L_i$ should be conjugate in a specific way:
\begin{equation}
\label{eqn:all-meridians-same}
\begin{aligned}
 g_D g_A^{-1} &= \alpha g_C g_B^{-1} \alpha^{-1}		& g_D \beta^{-m} g_A^{-1} \alpha^{m} &= \alpha g_C \beta^{-m} g_B^{-1} \alpha^{m} \alpha^{-1}\\
 g_A g_B^{-1} g_C g_A^{-1} &= \alpha^m g_C \beta^{-m} g_B^{-1}	\hspace{0.2in}	& g_D \beta^{-m} g_B^{-1} \alpha^{m} g_C g_D^{-1} &=  g_C g_B^{-1}
\end{aligned}
\end{equation}

Finally, the edges $u_4 u_2$, $v_1 v_3$, $v_4' v_2'$ and $u_3' u_1'$ are all identified. We must enforce that this artificial edge is non-singular in the glued up manifold:
\begin{align}
\label{eqn:artificial-edge} \alpha g_B \beta g_A^{-1} &= 1 	&  \alpha g_C \beta g_D^{-1} &= 1
\end{align}

Next, we convert these equations into equations about the (holonomy representation of) words in the fundamental group. The following is a generating set for $\pi_1\ N \setminus \Sigma$.
\begin{align*}
\mu &= g_D g_A^{-1} & \widetilde \beta &= g_A \beta g_A^{-1} & t &= g_B g_A^{-1}
\end{align*}
Then a quick computation gives that 
$$g_C g_A^{-1} = g_C g_B^{-1} g_B g_A^{-1} =  \alpha^{-1} \mu \alpha t$$

So we can easily translate any of the above equations into these new letters by inserting $g_A^{-1} g_A$ and using the above relations. 
We now transform the equations and reduce. First, we use either of (the equivalent) Equations~(\ref{eqn:artificial-edge}) to eliminate $t$:
\begin{align*}
\alpha t \widetilde \beta &= 1  & \implies& & t &= (\widetilde \beta \alpha)^{-1}
\end{align*}
Equations~(\ref{eqn:order2-fiber}) become:
\begin{align}
\label{will-give-fiber-length}
(\widetilde \beta \alpha)^2 &= \widetilde \beta^m \mu^{-1} \alpha^{-1} \mu \alpha 		& (\alpha \widetilde \beta)^2 &= \alpha^m \mu \widetilde \beta^{-1} \mu^{-1}\widetilde \beta 
\end{align}

Equations~(\ref{eqn:all-meridians-same}) now reduce to one equation. The equations on the first line are implied by the above, while the two equations of the second line turn out to be equivalent giving:
\begin{align} \label{m-eqn}
\mu^2 &= \alpha^{-m} \widetilde \beta^{m}
\end{align}
We demonstrate the derivation of this equation starting from the bottom right equation of~(\ref{eqn:all-meridians-same}):
\begin{align*}
g_D \beta^{-m} g_B^{-1} \alpha^{m} g_C g_D^{-1} &=  g_C g_B^{-1} \\
\mu \widetilde \beta^{-m} t^{-1} \alpha^m (\alpha^{-1} \mu \alpha t) \mu^{-1} &= (\alpha^{-1} \mu \alpha t) t^{-1} \\
 \widetilde \beta^{-m+1}  \alpha^m  \mu \widetilde \beta^{-1} &= \mu^{-1} \alpha^{-1} \mu \alpha \mu\\
  \widetilde \beta^{-m+1}  \alpha^m  \mu \widetilde \beta^{-1} &=  \alpha^{-1} \mu \alpha & \text{(Lemma~\ref{lem:inf-commute})}\\
 \widetilde \beta^{-m} \alpha^m \mu &= \beta^{-1} \alpha^{-1}\mu \alpha \beta \\
  \widetilde \beta^{-m} \alpha^m \mu &= \mu^{-1}. &\text{(Lemma~\ref{lem:kHkinv})}
\end{align*}

Next, we use Equations~\ref{will-give-fiber-length} and~\ref{m-eqn} to derive the relationship between $\dot \theta_\alpha$, $\dot \theta_\beta$, and $\dot \theta_{\beta \alpha}$, where $\dot \theta_{\beta \alpha}$ refers to the rate of change of the rotation angle for $\wt \beta \alpha$. Substituting from Equation~\ref{m-eqn} into the lefthand equation of~\ref{will-give-fiber-length} we have:
\begin{align*}
\widetilde \beta^{-m} (\widetilde \beta \alpha)^2 &= \sqrt{\alpha^m \widetilde \beta^{-m}}\ \alpha^{-1} \sqrt{\alpha^{-m} \widetilde \beta^m } \ \alpha \\
&= \sqrt{ \widetilde \beta^{-m}}  \sqrt{\alpha^m} \ \alpha^{-1} \sqrt{\alpha^{-m}} \sqrt{ \widetilde \beta^m } \ \alpha \\
&=  \sqrt{ \widetilde \beta^{-m}} \ \left(\alpha^{-1}\sqrt{ \widetilde \beta^m } \ \alpha \right)\\
\implies \ \ (\widetilde \beta \alpha)^2 &=  \sqrt{ \widetilde \beta^{m}} \ \left(\alpha^{-1}\sqrt{ \widetilde \beta^m } \ \alpha \right)
\end{align*}
where we note that if $k = 1+a\sigma$ is an infinitesimal, then $\sqrt{k} = 1 + \frac{1}{2}a \sigma$ is well-defined.
Now, $\widetilde \beta^m, \alpha^{-1} \widetilde \beta^m \alpha$ and $(\widetilde \beta \alpha)^2$ are infinitesimal rotations of $\HH^2$ about fixed points $\widetilde q, A^{-1} \widetilde q$ and $r$ respectively, where recall $A$ is the finite part of $\alpha$. Note that $r$ is the midpoint of the segment connecting $\widetilde q$ and $A^{-1} \widetilde q$ (see Figure~\ref{2mmglueing2}). Hence, $\sqrt{ \widetilde \beta^{m}} \ \left(\alpha^{-1}\sqrt{ \widetilde \beta^m } \ \alpha \right)$ is an infinitesimal rotation centered at $r$ with infinitesimal rotation amount $\cosh d(q,r) \ m \dot \theta_\beta$. So, the equation is satisfied if and only if $$2 \dot \theta_{\beta \alpha} = m \cosh d(q,r) \  \dot \theta_\beta.$$

Similarly, the righthand equation of~\ref{will-give-fiber-length} gives:
$$2 \dot \theta_{\alpha \beta} = m \cosh d(q,r) \  \dot \theta_\alpha$$
where $\dot \theta_{\alpha \beta}$ denotes the rate of change of rotation angle for $\alpha \widetilde \beta$.
From~(\ref{m-eqn}) we also get the relationship between the infinitesimal cone angle $\phi$ of $\mu$ and the rate of change of rotational parts of the generators. Note that, by Lemma~\ref{lem:HP3-action}, the infinitesimal cone angle $\phi$ is given by the infinitesimal translational part of $\mu$ along the projection of the singular axis to $\HH^2$: 
\begin{align*}
 2 \phi &= -2 m \frac{\sinh d(q,r)}{\sin \frac{\pi}{4}} \dot \theta_{\beta} \\
 &= -2 \sqrt{2} m \dot \theta_\beta  \sqrt{ \cot^2 \frac{\pi}{m} -1}
 \end{align*}

We have assumed that $\dot \theta_\beta > 0$ throughout. It is possible to interpret the case $\dot \theta_\beta < 0$ as building the same manifolds but with the fiber direction reversed (the building blocks in Figure~\ref{2mm-FPFQ} turn upside down). In this case, the infinitesimal rotational part of $\mu$ also changes sign (as in the above equation). However, the infinitesimal cone angle is always negative. For, if the orientation of the fibers is reversed, $\mu$ will encircle the singular locus in the negative direction; one should then use $\mu^{-1}$ to measure the infinitesimal cone angle.

\subsubsection{Placement of the $L_i$}

The last step in the construction of our singular $\HP$ structure on $N$, is to determine the exact placement of the lines $L_i$ on the boundary of the two building blocks for our structure. When glued together, all $L_i$ are matched up by the glueing maps, and the result must be a closed geodesic in the $\HP$ structure (with no corners). We can deduce the correct geodesic candidate for say $L_2$, by beginning at a point near $L_2$ and following once around the singular locus. The resulting group element $$g_B \beta^{-1} \beta^m g_D^{-1} \alpha^{-1} = \alpha^{-1} \widetilde \beta^{-1} \alpha \widetilde \beta$$ is a longitude of the singular locus, and is represented by a translation in $\HP$ whose geodesic axis must be $L_2$. This determines the other $L_i$ according to the equations in the previous section. The $H_i$ are also determined, once a choice of initial glueing map (say $g_A$) is made. We must check that the axes $L_i$ are actually arranged as in Figure~\ref{2mmglueing1}.

\begin{Lemma}
\label{lem:placement-Li}
\begin{enumerate}
\item
$L_2$ lies strictly in between $L_1$ and $L_3$ in the degenerate plane $\pi^{-1} e$ and does not intersect either.
\item 
$L_4$ lies strictly in between $L_1$ and $L_2$.
\end{enumerate}
\end{Lemma}
\begin{proof}
By symmetry ($\theta_\alpha = \theta_\beta$ and $\dot \theta_\alpha = \dot \theta_\beta$), we may assume that $\widetilde \beta = R \alpha R^{-1}$, where $R$ is rotation by $\pi$, preserving a degenerate plane containing $L_2$ as in Figure~\ref{proofdiagram2mm}. Then, with notation as in the figure, note that
\begin{align*}
w_1 &= \alpha x_2 \\
y_1 &= \widetilde \beta z_2 = R \alpha R^{-1} z_2 = R \alpha x_2 = R w_1
\end{align*}
\begin{figure}[h]
{\centering

\def\svgwidth{2.3in}
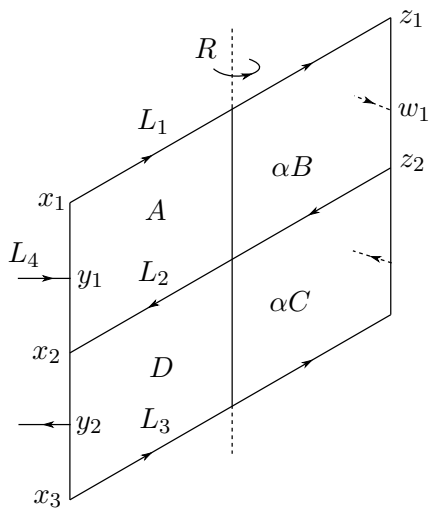

\caption{diagram for proof of Lemma~\ref{lem:placement-Li} \label{proofdiagram2mm}}
}
\end{figure}

The segment ${z_2 w_1}$ therefore has the same signed length as ${x_2 y_1}$. So $y_1 x_1 = \beta z_2 w_1$ then has the same signed length as $x_2 y_1$. So all four segments $x_3 y_2$, $y_2 x_2$, $x_2 y_1$ and $y_1 x_1$ have the same signed length. They must all be nontrivial, and positive since $x_1 = \alpha^m x_3 > x_3$.
\end{proof}

\subsection{Regeneration}

In the previous section, we showed that there exists an $\HP$-structure on $N$ with prescribed infinitesimal cone angle $\phi < 0$ along $\Sigma$. Let $M = N \setminus \Sigma$. The following Proposition shows that the cohomology condition of Theorem~\ref{thm:main-regen} is satisfied by our $\HP$ structure. This will complete the proof of Theorem~\ref{thm:regen-2mm}.

\begin{Proposition}\label{prop:2mm-H1}
The $\HP$ structure on $N$ with infinitesimal cone angle $\phi < 0$ along $\Sigma$ defined in the previous section is locally (and infinitesimally) rigid, if the infinitesimal cone angle $\phi$ is fixed. In particular, if $\rho = (\rho_0, z): \pi_1 M \to \PSL(2,\RR) \ltimes \mathfrak{sl}(2,\RR)$ is the holonomy representation, then $$H^1(\pi_1M, \mathfrak{sl}(2,\RR)_{Ad \rho_0} ) = \RR$$ is spanned by $z$.
\end{Proposition}

\begin{proof}
All quantities in the construction of the previous section are determined by the choice of one parameter, say, $\dot \theta_\beta$, which is related to $\phi$ by the linear equation
\begin{align*}
 \phi  &= - \sqrt{2} m \dot \theta_\beta  \sqrt{ \cot^2 \frac{\pi}{m} -1}.
 \end{align*}
when $\dot \theta_\beta > 0$.
So the first assertion will follow if we argue that any small change in the geometry corresponds to an $\HP$ structure which can still be constructed in the same way. Any $\HP$ structure on $N$ with infinitesimal cone singularity along $\Sigma$ gives rise via projection to a transversely hyperbolic foliation on $N$ (which extends smoothly over the $\Sigma$). By a theorem from Thurston's notes (Thm 4.9 \cite{Thurston}, see Section~\ref{hyperbolic-foliations}), this transversely hyperbolic foliation must be the same one we started with, namely that arising from projection of $N$ onto $S$, because the geometry of $S$ is rigid. So the $\PSL(2,\RR)$ part of the holonomy representation is rigid. Then, we may decompose our perturbed $\HP$ structure into two pieces as in Figure~\ref{2mmglueing1} by cutting the surface $S$ into two pieces and taking the inverse image under the projection as in the beginning of this section. This proves the first assertion about local rigidity and infinitesimal rigidity rel $\phi$.

The cohomology assertion follows because $\HP$ structures are locally parameterized by their holonomy representations which, since $\rho_0$ is rigid, are described by the group $H^1(\pi_1M, \mathfrak{sl}(2,\RR)_{Ad \rho_0} )$. Locally, classes $[z]$ are in $1-1$ correspondence with infinitesimal cone angles $\phi$ (though, note that $[z]$ and $[-z]$ correspond to the same infinitesimal cone angle).
\end{proof}

A generalization of Theorem~\ref{thm:regen-2mm} and the methods of this section to the case of higher genus $S$ may be possible. However, it is important to note the following obstacle. The $\PSL(2,\RR)$ representation variety will not be smooth in such a case: a representation corresponding to a collapse of the three manifold onto a surface $S$ lies on the intersection of two irreducible components, the extra one corresponding to deforming the geometry of $S$ (as the Teichmuller space of $S$ is non-trivial in this case). In the final section, we give an example displaying this phenomenon and some interesting consequences of it.

%
%

\section{An interesting flexibility phenomenon}
\label{s:brings}
Here we will construct examples of transitioning structures for which the underlying hyperbolic and $\AdS$ structures collapse onto a puncture torus. The underlying manifold is the Borromean rings complement $M$ with one boundary torus required to be a parabolic cusp and the other two boundary tori filled in with cone/tachyon singularities. We will see that the $\SO(2,1)$ representation variety is not smooth at the locus of degenerated structures, so the existence of a representation-level transition is not automatic. As a result of the non-smoothness we observe that a transitional $\HP$ structure on $M$ can be deformed to nearby $\HP$ structures that \emph{do not} regenerate to hyperbolic structures. Interestingly,  
these nearby $\HP$ structures do regenerate to $\AdS$ structures. The geometry can be constructed using ideal tetrahedra and the methods of \cite{Danciger-11} (in fact, $M$ is the union of eight tetrahedra). However, for brevity, we observe this phenomenon only at the level of representations.

\subsubsection{Representation variety}
\begin{figure}[h]
\centering
\includegraphics[width=2.0in]{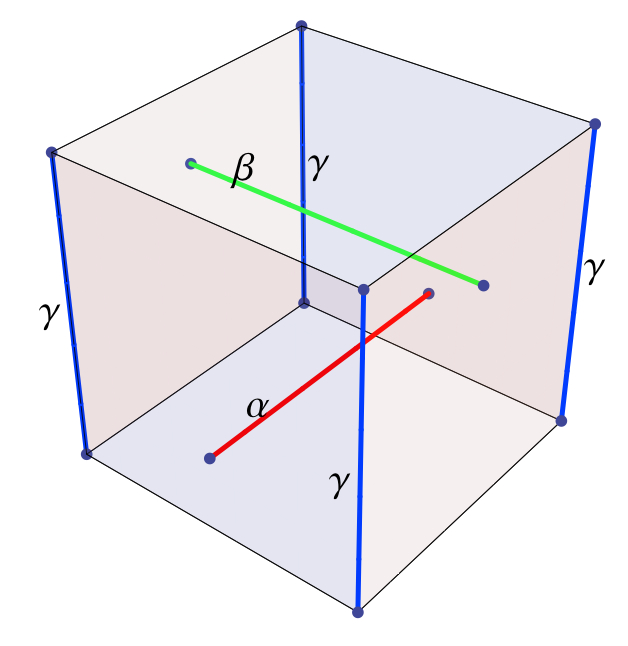}
\caption[The Borromean rings complement.]{We remove the three curves $\alpha, \beta, \gamma$ shown in the diagram from the three-torus $T^3$ (opposite sides of the cube are identified). The resulting manifold $M$ is homeomorphic to the complement of the Borromean rings in $S^3$.}
\label{brings}
\end{figure}

Consider the three-torus $T^3$ defined by identifying opposite faces of a cube. Now, define $M^3 = T^3-\{\alpha,\beta,\gamma\}$, where $\alpha, \beta, \gamma$ are disjoint curves freely homotopic to the generators $a,b,c$ of $\pi_1 T^3$ as shown in Figure~\ref{brings}.
 Then $M$ is homeomorphic to the complement of the Borromean rings in $S^3$ (this is stated in \cite{Hodgson-86}). A presentation for $\pi_1 M$ is given by: $$\pi_1 M = \langle a,b,c : [[a,b],c] = [[c,b^{-1}],a] = 1 \rangle.$$
We study the representation variety $\mathscr R_{par}(M)$ of representations $\rho: \pi_1 M \rightarrow \PSL(2, \mathbb R)$ up to conjugacy such that $\rho[a,b]$ is parabolic (and so $\rho(c)$ is parabolic with the same fixed point). These representations correspond to transversely hyperbolic foliations which are ``cusped" at one boundary component and have Dehn surgery type singularities at the other two boundary components (see \cite{Hodgson-86}).

Let $T$ denote the punctured torus, with $\pi_1 T = \langle a, b\rangle$. Then $\pi_1 T \hookrightarrow \pi_1 M$, so that $\mathscr R_{par}(M) \rightarrow \mathscr R_{par}(T)$. The representations of $\mathscr R_{par}(T)$ correspond to hyperbolic punctured tori (with a cusp at the puncture). A representation $\rho : \pi_1 T \rightarrow \PSL_2 \mathbb R$ satisfies the parabolic condition if and only if $\rho(a), \rho(b)$ are hyperbolic elements with $$\sinh\frac{l(a)}{2}\sinh\frac{l(b)}{2}\sin\phi = 1 $$ where $l(a), l(b)$ are the translation lengths of $\rho(a), \rho(b)$ respectively and $\phi$ is the angle between the axes. To lift such a representation to a representation of $\pi_1 M$, we must assign $\rho(c)$ so that the relations of $\pi_1 M$ are satisfied. Since $\rho(c)$ must commute with the parabolic $\rho[a,b]$, $\rho(c)$ is parabolic with the same fixed point. Let $x$ denote the amount of parabolic translation of $\rho(c)$ relative to $\rho[a,b]$, so if $\rho[a,b] = \left(\begin{smallmatrix} -1 & 1\\ 0 & -1\end{smallmatrix}\right)$, then $\rho(c) = \left(\begin{smallmatrix} -1 & x\\ 0 & -1\end{smallmatrix}\right)$. It turns out (by a nice geometric argument) that there are exactly two solutions for $x$: 
\begin{equation*}
x = 0 \ \ \ \text{ or } \ \ x = \frac{1}{2}\text{sech}\frac{l(a)}{2}\text{sech}\frac{l(b)}{2}\cot \phi.
\end{equation*} 
This describes the representation variety $\mathscr R_{par}$ rather explicitly as the union of two irreducible two-dimensional components $\mathscr R_T$ and $\mathscr R_R$. The first component $\mathscr R_T$ (`T' for Teichmuller) consists of the obvious representations with $\rho(c) = 1$ and $\rho(a),\rho(b)$ generating a hyperbolic punctured torus group. The associated transversely hyperbolic foliations are products (with two fillable singularities at $\alpha$ and $\beta$). The second component $\mathscr R_R$ (`R' for regenerate) describes transversely hyperbolic foliations with more interesting structure. This component, in fact its complexification, is the relevant one for regenerating hyperbolic structures. Note that $\mathscr R_T$ and $\mathscr R_R$ meet (transversely) exactly at the locus of ``rectangular" punctured tori ($\cot \phi = 0$).
\begin{Remark}
If we identify $\mathscr R_T$ with the Teichmuller space $\mathcal T_{1,1}$ of the punctured torus, then the singular set of $\mathscr R_{par}$, given by $\mathscr R_T \cap \mathscr R_R$, is exactly the \emph{line of minima} for the curves $a$ and $b$. In other words, $\mathscr R_T \cap \mathscr R_R$ consists of the representations in $\mathscr R_T$ where there is a relation between the differentials $dl(a)$ and $dl(b)$. The relevance of such a relation in the context of regeneration questions is discussed in Section 3.17 of \cite{Hodgson-86}. 
\end{Remark}


\subsubsection{Regenerating 3D structures}
Fix a particular rectangular punctured torus $\rho_0 : \pi_1 T \rightarrow \PSL(2,\RR) \mathbb R$, and lift $\rho_0$ to $\pi_1 M$ by setting $\rho_0(c) = 1$ (this is the only possible lift). Let $v$ be a tangent vector at $\rho_0$, tangent to the component $\mathscr R_R$ but transverse to $\mathscr R_T$. For suitably chosen $v$, the representation $\rho_0 + \sigma v : \pi_1 M \rightarrow \PSL_2(\mathbb R + \mathbb R \sigma)$ is the holonomy of a robust $\HP$ structure (which can be constructed from eight tetrahedra). Now, as the variety $\mathscr R_R$ is smooth, the complexified variety $\mathscr R_R^{\mathbb C}$ is smooth at $\rho_0$. Thus the Zariski tangent vector $iv$ is tangent to a path $\rho_t : \pi_1 M \rightarrow \PSL(2, \mathbb C)$ which is compatible to first order with $\rho_0 + \sigma v$. By Proposition~\ref{regen}, the $\HP$ structure regenerates to a path of hyperbolic structures with holonomy $\rho_t$ (or alternatively, this path of hyperbolic structures can be constructed directly using tetrahedra). 
Similarly, the variety $\mathscr R_R^{\mathbb R + \mathbb R \tau}$ is smooth at $\rho_0$ yielding a path of holonomies $\rho_t : \pi_1 M \rightarrow \PSL(2,\Rtau)$ with $\rho'_0 = \tau v$ so that Proposition~\ref{regen} then produces a regeneration to $\AdS$ structures with holonomy $\rho_t$. Thus, our $\HP$ structure is transitional.
Actually, in the $\AdS$ case, the representations can be constructed directly. Let $\sigma_t : \pi_1 M \rightarrow \PSL(2, \mathbb R)$ be a path with $\sigma'_0 = v$. Then, a path $\rho_t$ of $\PSL(2, \Rtau)$ representations with $\rho'_0 = \tau v$ is defined by $$\rho_t = \frac{1+\tau}{2} \sigma_t + \frac{1-\tau}{2} \sigma_{-t}.$$ 

\subsubsection{An interesting flexibility phenomenon}
\begin{figure}[h]
{\centering

\def\svgwidth{2.5in}
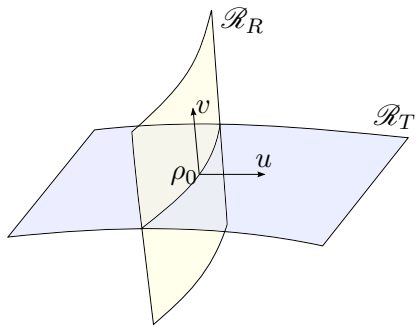

}

\caption[A schematic picture of the $\PSL(2,\mathbb R)$ representation variety $\mathscr R_{par}(M)$.]{A schematic picture of the $\PSL(2,\mathbb R)$ representation variety $\mathscr R_{par}(M)$. The variety is the union of two irreducible two-dimensional components which meet at the locus of rectangular punctured torus representations (with $c = 1$). We let $\rho_0$ be one such representation, with $v$ tangent to one component, and $u$ tangent to the other.}
\label{tangent-vectors}
\end{figure}

The transitional $\HP$ structure from the previous sub-section, with holonomy $\rho_0 + \sigma v$, can be deformed in an interesting way. By Proposition~\ref{prop:deform-withboundary}, nearby $\HP$ structures are determined by nearby holonomy representations. We consider a deformation of the form $$\rho_0 + \sigma (v + \epsilon u)$$ where $\epsilon u$ is a small tangent vector at $\rho_0$, tangent to the component $\mathscr R_T$ and transverse to $\mathscr R_R$ (see Figure~\ref{tangent-vectors}).
 Now, $\mathscr R_{par}^{\mathbb C}(M)$ is the union of its irreducible components $\mathscr R_T^{\mathbb C}$ and $\mathscr R_R^{\mathbb C}$ (locally at $\rho_0$). So, as $u$ and $v$ are tangent to different components of $\mathscr R_{par}(M)$, any Zariski tangent vector of the form $w+i(v + \epsilon u)$, for $w$ real, is not integrable. Thus, the deformed $\HP$ structure \emph{does not} regenerate to hyperbolic structures. 
 %
However, it does regenerate to $\AdS$ structures. To see this, consider paths $\sigma_t$ and $\mu_t$ with derivatives $2v$ and $2\epsilon u$ respectively at $t=0$. Then, $$\rho_t = \frac{1+\tau}{2} \sigma_{t} + \frac{1-\tau}{2} \mu_{-t}$$ gives a family of $\PSL(2,\Rtau)$ representations with $\rho'_0 = (v-\epsilon u) + \tau(v+\epsilon u)$. Proposition~\ref{regen} now implies that the deformed $\HP$ structure regenerates to $\AdS$ structures.

\bibliographystyle{amsalpha}
\small
\bibliography{mybib}

\newcommand{\etalchar}[1]{$^{#1}$}
\providecommand{\bysame}{\leavevmode\hbox to3em{\hrulefill}\thinspace}
\providecommand{\MR}{\relax\ifhmode\unskip\space\fi MR }
\providecommand{\MRhref}[2]{%
  \href{http://www.ams.org/mathscinet-getitem?mr=#1}{#2}
}
\providecommand{\href}[2]{#2}
\begin{thebibliography}{BBD{\etalchar{+}}12}

\bibitem[ABB{\etalchar{+}}07]{Andersson-07}
Lars Andersson, Thierry Barbot, Riccardo Benedetti, Francesco Bonsante,
  William~M. Goldman, Francois Labourie, Kevin~P. Scannell, and Jean-Marc
  Schlenker, \emph{Notes on a paper of {M}ess}, Geometriae Dedicata
  \textbf{126} (2007), no.~1, 47--70.

\bibitem[BB09]{Benedetti-09}
Riccardo Benedetti and Francesco Bonsante, \emph{Canonical wick rotations in
  3-dimensional gravity}, Memoirs of the American Mathematical Society, 2009.

\bibitem[BBD{\etalchar{+}}12]{BBDGGKKSZ-12}
Thierry Barbot, Francesco Bonsante, Jeffrey Danciger, William~M. Goldman,
  Fran\c cois Gu\'eritaud, Fanny Kassel, Kirill Krasnov, Jean-Marc Schlenker,
  and Abdelghani Zeghib, \emph{Some questions on anti de-sitter geometry},
  arXiv:1205.6103 \url{http://arxiv.org/abs/1205.6103} (2012).

\bibitem[BBS09]{Barbot-09}
Thierry Barbot, Francesco Bonsante, and Jean-Marc Schlenker, \emph{Collisions
  of particles in locally {A}d{S} spacetimes}, preprint (2009).

\bibitem[BEE96]{Beem-book}
John~K. Beem, Paul~E. Ehrlich, and Kevin~L. Easley, \emph{Global lorentzian
  geometry, 2nd ed.}, Marcel Dekker, New York, 1996.

\bibitem[Ben08]{Benoist-04}
Yves Benoist, \emph{A survey on divisible convex sets}, Adv. Lect. Math
  \textbf{6} (2008), 1--18.

\bibitem[Ber60]{Bers-60}
Lipman Bers, \emph{Simultaneous uniformization}, Bull. Amer. Math. Soc.
  \textbf{66} (1960), 94--97.

\bibitem[BLP05]{Boileau-05}
Michel Boileau, Bernhard Leeb, and Joan Porti, \emph{Geometrization of
  3-dimensional orbifolds}, Annals of Math \textbf{162} (2005), no.~1,
  195--250.

\bibitem[Bro07]{Bromberg-07}
Ken Bromberg, \emph{Projective structures with degenerate holonomy and the bers
  density conjecture}, Annals of Math \textbf{166} (2007), 77--93.

\bibitem[CG97]{Schoi-97}
Suhyoung Choi and William~M. Goldman, \emph{The classification of real
  projective structures on surfaces}, Bulletin of the American Math Society
  \textbf{34} (1997), 161--171.

\bibitem[CHK00]{Cooper-00}
Daryl Cooper, Craig Hodgson, and Steven Kerckhoff, \emph{Three-dimensional
  orbifolds and cone manifolds}, MSJ Memoirs, vol.~5, Mathematical Society of
  Japan, Tokyo, 2000.

\bibitem[CLT12]{Cooper-12}
Daryl Cooper, Darren~D. Long, and Stephan Tillmann, \emph{On convex projective
  manifolds and cusps}, preprint (2012).

\bibitem[CM12]{Crampon-12}
Micka{\"e}l Crampon and Ludovic Marquis, \emph{Finitude g\' eom\' etrique en
  g\' eom\' etrie de hilbert}, arXiv: 1202.5442v1 (2012).

\bibitem[Dan11]{Danciger-11}
Jeffrey Danciger, \emph{Geometric transitions: from hyperbolic to {AdS}
  geometry}, Ph.D. thesis, Stanford University, 2011.

\bibitem[Ehr36]{Ehresmann-36}
Charles Ehresmann, \emph{Sur les espaces localement homogenes}, L'Enseignement
  MathŽmatique \textbf{35} (1936), 317--333.

\bibitem[Gol90]{Goldman-90}
William~M. Goldman, \emph{Convex real projective structures on compact
  surfaces}, Journal of Differential Geometry \textbf{31} (1990), 791--845.

\bibitem[Gol10]{Goldman-10}
\bysame, \emph{Locally homogeneous geometric manifolds}, Proceedings of the
  International Congress of Mathematicians, Hyderabad, India (2010).

\bibitem[HK98]{Hodgson-98}
Craig Hodgson and Steven Kerckhoff, \emph{Rigidity of hyperbolic cone manifolds
  and hyperbolic dehn surgery}, Journal of Differential Geometry \textbf{48}
  (1998), 1--60.

\bibitem[HK05]{Hodgson-05}
\bysame, \emph{Universal bounds for hyperbolic dehn surgery}, Annals of
  Mathematics \textbf{162} (2005), no.~1, 367--421.

\bibitem[Hod86]{Hodgson-86}
Craig Hodgson, \emph{Degeneration and regeneration of hyperbolic structures on
  three-manifolds}, Ph.D. thesis, Princeton University, 1986.

\bibitem[HPS01]{Huesener-01}
Michel Huesener, Joan Porti, and Eva Su\'arez, \emph{Regenerating singular
  hyperbolic structures from sol}, Journal of Differential Geometry \textbf{59}
  (2001), 439--478.

\bibitem[Mes07]{Mess-07}
Geoffrey Mess, \emph{Lorentz spacetimes of constant curvature}, Geometriae
  Dedicata \textbf{126} (2007), no.~1, 3--45.

\bibitem[Por98]{Porti-98}
Joan Porti, \emph{Regenerating hyperbolic and spherical cone structures from
  euclidean ones}, Topology \textbf{37} (1998), no.~2, 365--392.

\bibitem[Por02]{Porti-02}
\bysame, \emph{Regenerating hyperbolic cone structures from nil}, Geometry and
  Topology \textbf{6} (2002), 815--852.

\bibitem[Por10]{Porti-10}
\bysame, \emph{Regenerating hyperbolic cone 3-manifolds from dimension 2},
  preprint (2010).

\bibitem[Rat94]{Ratcliffe-book}
John~G. Ratcliffe, \emph{Foundations of hyperbolic manifolds}, Springer-Verlag,
  1994.

\bibitem[Suo03]{Suoto-03}
Juan Suoto, \emph{Hyperbolic cone-manifolds with large cone-angles}, Geometry
  and Topology \textbf{7} (2003), 789--797.

\bibitem[Thu80]{Thurston}
William~P. Thurston, \emph{The geometry and topology of three manifolds},
  \url{www.msri.org/publications/books.gt3m}, 1980.

\bibitem[Thu97]{Thurston-book}
\bysame, \emph{Three-dimensional geometry and topology}, Princeton University
  Press, 1997.

\end{thebibliography}

\end{document}